\documentclass[12pt]{article}
\usepackage{setspace}
\usepackage{amsmath,amssymb,amsthm}
\usepackage{mathrsfs}
\usepackage{geometry}
\usepackage{enumitem}
\usepackage{environ}
\usepackage{setspace}
\usepackage{titlesec}
\usepackage[all,cmtip]{xy}
\usepackage{tikz}
\usepackage{graphicx}
\usetikzlibrary{matrix,arrows}

\NewEnviron{prf}[1][]{\begin{proof}[\bf #1Proof]\BODY\end{proof}}{}
\NewEnviron{slt}[1][]{\begin{proof}[\bf #1	解]\BODY\end{proof}}{}
\newtheorem{definition}{Definition}[section]
\newtheorem{theorem}[definition]{Theorem}
\newtheorem{lemma}[definition]{Lemma}
\newtheorem{proposition}[definition]{Proposition}
\newtheorem{remark}[definition]{Remark}

\newtheorem{corollary}[definition]{Corollary}

\makeatletter
\newcommand*\bigcdot{\mathpalette\bigcdot@{.5}}
\newcommand*\bigcdot@[2]{\mathbin{\vcenter{\hbox{\scalebox{#2}{$\m@th#1\bullet$}}}}}
\makeatother
  
\setlist[enumerate,1]{label=(\roman*)}
\setlist[enumerate,2]{label=(\arabic*)}

\title{\textbf{A regularization theorem for bounded-degree self-maps}}
\author{She YANG}
\date{}

\begin{document}
\begin{spacing}{1.25}

\maketitle

\begin{abstract}
Let $K$ be an algebraically closed field of arbitrary characteristic and let $X$ be an irreducible projective variety over $K$. Let $G\subseteq\text{Bir}(X)$ be a bounded-degree subgroup. We prove that there exists an irreducible projective variety $Y$ birational to $X$, such that every element of $G$ becomes an automorphism of $Y$ after the birational transformation. If $K=\mathbb{C}$, this result is stated in [Can14, Theorem 2.5] and the proof backs to [HZ96, Section 5]. The proof in [HZ96] is not purely algebraic. Inheriting the methods in [HZ96], we give a purely algebraic proof of this statement in arbitrary characteristic. We will also discuss a corollary of this result which is useful in arithmetic dynamics.
\end{abstract}

\section{Introduction}

In this note, we fix an algebraically closed field $K$ of arbitrary characteristic. A variety is a reduced separated finite type $K$-scheme. The fiber products are over $K$ if not specified.

Let $X$ be an irreducible projective variety. We denote $\text{Bir}(X)$ as the group of birational self-maps of $X$. Let $G\subseteq\text{Bir}(X)$ be a subgroup. It is natural to ask that whether there exists an irreducible projective variety $Y$ birational to $X$, such that every element of $G$ becomes an automorphism of $Y$ after the birational transformation. The Weil's regularization theorem $\cite[\text{Theorem}\ 2.5]{Can14}$ says that if $K=\mathbb{C}$ and $G$ is a \emph{bounded-degree} subgroup, then the model $Y$ exists. We will prove this statement for arbitrary algebraically closed field $K$.

We explain the meaning of bounded-degree. Let $X$ be an irreducible projective variety and let $L\in\text{Pic}(X)$ be an ample line bundle. For a self-map $f\in\text{Bir}(X)$, we denote $\Gamma_{f}$ as the closed graph of $f$, which is an irreducible closed subvariety of $X\times X$. Let $L_{0}=\text{pr}_{1}^{*}(L)+\text{pr}_{2}^{*}(L)$ be an ample line bundle on $X\times X$ where $\text{pr}_{1},\text{pr}_{2}:X\times X\rightarrow X$ are the two projections. Then we let the \emph{degree} of $f$ with respect to $L$ be the intersection number $(L_{0}^{\text{dim}(X)}\cdot\Gamma_{f})$ and denote as $\text{deg}_{L}(f)$.

\begin{definition}
Let $X$ be an irreducible projective variety. A subgroup $G\subseteq\mathrm{Bir}(X)$ is of \emph{bounded-degree} if $\{\mathrm{deg}_{L}(f)|\ f\in G\}$ is bounded for some (and hence every) ample line bundle $L\in\mathrm{Pic}(X)$.
\end{definition}

Our main theorem is as follows.

\begin{theorem}
Let $X$ be an irreducible projective variety and let $G\subseteq\mathrm{Bir}(X)$ be a bounded-degree subgroup. Then there exists an irreducible projective variety $Y$ birational to $X$, such that every element of $G$ becomes an automorphism of $Y$ after the birational transformation.
\end{theorem}

There is a corollary of this result which is useful in arithmetic dynamics. It says that bounded-degree self-maps are regularizable. We will state the definition of bounded-degree self-maps and prove this corollary in the last section.

\begin{corollary}
Let $X$ be an irreducible projective variety and let $f$ be a bounded-degree self-map of $X$. Then there exists an irreducible projective variety $Y$ birational to $X$, such that $f$ becomes a bounded-degree automorphism of $Y$ after the birational transformation.
\end{corollary}

To prove Theorem 1.2, we mainly follow the route in $\cite[\text{Section}\ 5]{HZ96}$ and there is a short summary of the proof in $\cite[\text{Subsection}\ 2.3]{Can14}$. The original approach in $\cite{HZ96}$ uses Chow variety as the moduli space while $\cite{Can14}$ uses Hilbert scheme. It turns out that in fact the argument can be carried out using each of them. We will give a detailed proof of Theorem 1.2 in Section 3 using Chow varieties, and then in Section 4 describe the changes in the details for another similar approach using Hilbert schemes.

Now we briefly summarize the idea of the proof of Theorem 1.2 using Chow varieties. Firstly, we will identify the birational self-maps in $G$ with their closed graphs. Since $G$ is of bounded-degree, the closed graphs $\{\Gamma_{f}\subseteq X\times X|\ f\in G\}$ will lie in finitely many open Chow varieties $C_{\text{dim}(X),d}(X\times X)$ (see Theorem 2.5 for the notion of open Chow variety). Now let $T$ be the closure of $G$ in $\bigsqcup\limits_{d}C_{\text{dim}(X),d}(X\times X)$, which is a quasi-projective variety. Let $S$ be the union of the irreducible components of maximal dimension in $T$, which is a closed subvariety of $T$. We will show that one can induce a rational group structure (see Definition 2.1) on $S$ by the group structure of $G$, and there is a natural rational group action of $S$ on $X$. Then the original Weil's regularization theorems about rational groups and rational group actions (see Theorem 2.4) can be applied, and hence the result follows.

~

At the end of the Introduction, we describe the structure of this note. We do some preparations in Section 2. We will introduce the two technical tools, i.e. the Weil's regularization theorems and the open Chow variety; and we will also propose several technical lemmas. Those lemmas can also be regarded as a warm-up before the more complicated proof of Theorem 1.2 in Section 3. Then in Section 4, we modify the argument in Section 3 to describe the proof of Theorem 1.2 using Hilbert schemes as the parameter space. Lastly, in Section 5, we will discuss about the bounded-degree self-maps in arithmetic dynamics and prove Corollary 1.3.

\clearpage

\textbf{Notation and Terminology.}

We fix an algebraically closed field $K$ of arbitrary characteristic. A variety is a reduced separated finite type $K$-scheme. The morphisms are $K$-morphisms, and the fiber products are over $K$ if not specified.

Let $X,Y$ be two varieties. A rational map $f:X\dashrightarrow Y$ is an equivalence class of $(U,f_{U})$ where $U\subseteq X$ is an open dense subset and $f_{U}:U\rightarrow Y$ is a morphism, as usual. A rational map $f$ has a largest domain of definition, which will be denoted as $\text{Dom}(f)$. If $f:X\dashrightarrow Y$ is a rational map, we will denote $\Gamma_{f}$ as the closed graph of $f$. It is a closed subvariety of $X\times Y$.

Let $X,Y,S$ be varieties and $p:X\rightarrow S$, $q:Y\rightarrow S$ be morphisms. We will say a rational map $f:X\dashrightarrow Y$ is an $S$-rational map if $q\circ f=p$ (notice that rational maps can always composite with morphisms in this way). Sometimes we will consider the closed graph $\Gamma_{f}\subseteq X\times_{S}Y$. We will make sure that there is no ambiguity.

Let $f:X\dashrightarrow Y$ be a rational map between two varieties. We will say that $f$ is \emph{dominant}, if
\begin{enumerate}
\item
$\overline{f(\text{Dom}(f))}=Y$, and

\item
the inverse image of open dense subsets in $Y$ are open dense in $X$.
\end{enumerate}

Some remarks should be made towards the definition of dominant rational maps. In $\cite{Dem70}$ and $\cite{Kra}$, the dominant rational map only has to satisfy condition (i) and in $\cite{Zai95}$, it only has to satisfy condition (ii). However, it is condition (ii) that guarantees one can compose dominant rational maps. So the condition (ii) cannot be abandoned. On the other hand, it is condition (i) which makes dominant rational maps look like ``dominant" as we can imagine. So we still want to reserve condition (i).

An equivalent definition is made in $\cite[\text{Remark}\ 29.49.13\ \text{(0BX9)}]{Stacks}$. Let $X^{0},Y^{0}$ be the finite set of generic points of the irreducible components of $X,Y$ respectively. Then condition (i) is equivalent to $Y^{0}\subseteq f(X^{0})$ and condition (ii) is equivalent to $f(X^{0})\subseteq Y^{0}$. So our requirement of dominant is equivalent to saying $f(X^{0})=Y^{0}$, just as in $\cite[\text{Remark}\ 29.49.13\ \text{(0BX9)}]{Stacks}$. From this point of view, one can see that if the number of irreducible components of $X,Y$ are equal (or more specifically, if $X=Y$), then condition (i) implies condition (ii). So we do not have to be overcautious with this concept when dealing with rational self-maps.

Let $f:X\dashrightarrow Y$ be a dominant rational map between two varieties. We say that $f$ is \emph{birational}, if there exists a dominant rational map $g:Y\dashrightarrow X$ such that $g\circ f=\text{id}_{X}$ and $f\circ g=\text{id}_{Y}$. Equivalently, $f:X\dashrightarrow Y$ is birational if and only if it can be represented by an isomorphism $f_{U}:U\stackrel{\sim}\rightarrow V$ where $U,V$ are open dense subsets of $X,Y$ respectively.

~

\textbf{Acknowledgement.} The author is very grateful to his advisor Junyi Xie who tells him this topic and constantly supports him through the work. He is also very grateful to Xiangyu Pan, whose undergraduate thesis provides a good reference about the Weil's regularization theorems and also contains a brief summary about $\cite[\text{Section}\ 5]{HZ96}$. He is very thankful to J\'anos Koll\'ar for verifying the statement of the open Chow variety (Theorem 2.5). He would also like to thank Chengyuan Yang for some useful discussions.

\section{Preliminaries}

In this section, we will do some preparations before proving Theorem 1.2 in Section 3. In Subsection 2.1, we will introduce the two main tools that will be used in the proof, i.e. the Weil's regularization theorems and the open Chow variety. In Subsection 2.2, we will prove some technical lemmas which will be used in Section 3 implicitly or explicitly. They can also be regarded as a warm-up before Section 3. And in Subsection 2.3, we will use the Weil's regularization theorem to give a quick proof of a baby case of Theorem 1.2, that is, the case when $G$ is a finite group.

\subsection{Two main tools}

We introduce the Weil's regularization theorems and the open Chow variety in this subsection.

~

In short, the Weil's regularization theorems say that one can regularize rational groups and rational group actions. So we will start with two definitions.

\begin{definition}
Let $X$ be a variety. Let $m:X\times X\dashrightarrow X$ be a rational map. We say that $(X,m)$ is a \emph{rational group}, if
\begin{enumerate}
\item
The diagram
\[
\xymatrix{X\times X\times X \ar@{.>}[r]^-{\mathrm{id}_{X}\times m} \ar@{.>}[d]_{m\times\mathrm{id}_{X}} & X\times X \ar@{.>}[d]^{m} \\ X\times X \ar@{.>}[r]^{m} & X}
\]
is commutative, and

\item
The rational maps
$$
\Phi:X\times X\dashrightarrow X\times X,(x,y)\mapsto(m(x,y),y)
$$
$$
\Psi:X\times X\dashrightarrow X\times X,(x,y)\mapsto(x,m(x,y))
$$
are birational.
\end{enumerate}

We will say that $X$ is a rational group if there is no ambiguity about $m$.
\end{definition}

\begin{definition}
Let $G$ be a group variety and let $X$ be an irreducible variety. A \emph{rational group action} of $G$ on $X$ is a rational map $\rho:G\times X\dashrightarrow X$ such that

\begin{enumerate}
\item
The diagram
\[
\xymatrix{G\times G\times X \ar@{.>}[r]^-{\mathrm{id}_{G}\times\rho} \ar[d]_{m_{G}\times\mathrm{id}_{X}} & G\times X \ar@{.>}[d]^{\rho} \\ G\times X \ar@{.>}[r]^{\rho} & X}
\]
is commutative, and

\item 
The rational map $\tilde{\rho}:G\times X\dashrightarrow G\times X,(g,x)\mapsto(g,\rho(g,x))$ is birational.
\end{enumerate}
\end{definition}

\begin{remark}
\begin{enumerate}
\item
Both $m$ and $\rho$ in the two definitions must be dominant. This is because of the condition (ii) in the definitions respectively. As a result, all of the maps in the two diagrams in condition (i) are dominant. So there is no problem about composition.

\item
In Definition 2.2(ii), we have $\mathrm{Dom}(\tilde{\rho})=\mathrm{Dom}(\rho)$.
\end{enumerate}
\end{remark}

Now we can state the Weil's regularization theorems.

\begin{theorem}(Weil's regularization theorems)
\begin{enumerate}
\item
Let $(X,m)$ be a rational group. Then there exist a group variety $G$ and a birational map $f:G\dashrightarrow X$, such that the diagram below is commutative.
\[
\xymatrix{G\times G \ar[r]^-{m_{G}} \ar@{.>}[d]_{f\times f} & G \ar@{.>}[d]^{f} \\ X\times X \ar@{.>}[r]^-{m} & X}
\]

\item
Let $G$ be a group variety and let $X$ be an irreducible variety. Let $\rho:G\times X\dashrightarrow X$ be a rational group action of $G$ on $X$. Then there exist an irreducible projective $G$-variety $Y$ and a birational map $f:X\dashrightarrow Y$, such that the diagram below is commutative (where $a$ is the group action of $G$ on $Y$).
\[
\xymatrix{G\times X \ar@{.>}[r]^-{\rho} \ar@{.>}[d]_{\mathrm{id}_{G}\times f} & X \ar@{.>}[d]^{f} \\ G\times Y \ar[r]^-{a} & Y}
\]
\end{enumerate}
\end{theorem}

The most original reference of this theorem is $\cite{Wei55}$, and one can consult $\cite{Zai95}$ for an approach in modern variety-theoretic language. For (i), this statement is $\cite[\text{Theorem}\ 3.7]{Zai95}$ and it can also be regarded as a special case of $\cite[5.1,\text{Theorem}\ 5]{BLR90}$. For (ii), the statement is a consequence of the results in $\cite{Bri22}$. More precisely, if we drop the projective requirement on $Y$, then this result is proved in $\cite[\text{Section}\ 6]{Bri22}$ (see also $\cite[\text{Theorem}\ 4.9]{Zai95}$ and $\cite[\text{Theorem}\ 1]{Kra}$) and the original assertion can be found at $\cite[\text{middle}\ \text{of}\ \text{p.}\ 11]{Bri22}$. By taking $\cite[\text{Theorem}\ 1,\text{Theorem}\ 2]{Bri22}$ into account, we may further require $Y$ to be projective.

Now we introduce the open Chow variety. The Chow variety is a moduli space which parametrizes cycles of given dimension and degree in a projective variety. It was firstly constructed in $\cite{CvdW37}$ and its functoriality was seriously studied in $\cite[\text{I}.3\text{--}4]{Kol96}$. Moreover, this concept turns out to be sensitive with characteristic, according to the statements in $\cite[\text{I}.4]{Kol96}$. In this note, we will use a version about the ``open" part of the Chow variety described in $\cite[3.7.1]{Kol23}$.

\begin{theorem}(Open Chow variety; a special case of $\cite[3.7.1]{Kol23}$)

Let $X$ be a projective variety and let $L$ be a very ample line bundle on $X$. Fix a nonnegative integer $k$ and a positive integer $d$. Then the functor from $\mathrm{(Seminormal}\ \mathrm{varieties)}^{op}$ to $\mathrm{(Sets)}$ as below can be represented by a seminormal quasi-projective variety $C_{k,d}(X)$ (the open Chow variety):
$$
S \mapsto \left\{
\begin{array}{@{}l@{}}
  F\subseteq X\times S\text{ is a closed}\\
  \text{subvariety of }X\times S\\
\end{array}
\;\left|\;\;
\begin{array}{@{}l@{}}
  \text{$\forall s\in S$, the fiber $F_s\subseteq X_{k(s)}$ has a}\\
  \text{geometrically reduced open dense} \\
  \text{subset, and $F_{s}$ is of pure dimension} \\
  \text{$k$ and degree $d$; and the morphism} \\
  \text{$F\rightarrow S$ is universally open.} \\
\end{array}
\right.
\right\}
$$
\end{theorem}

\begin{remark}
\begin{enumerate}
\item
The functoriality is given by ``$F\subseteq X\times S\mapsto F_{S',\text{red}}\subseteq X\times S'$" where $S'\rightarrow S$ is a given morphism of seminormal varieties, and the degree is calculated by intersecting with the line bundle $L$. One can verify that this is indeed a functor. Notice that on the contrary to the functor of Hilbert schemes, this is a variety-theoretic (or set-theoretic, and essentially cycle-theoretic) functor rather than a scheme-theoretic one.

\item
Seminormal is a requirement that stronger than reduced and weaker than normal. For references, one can consult $\cite[\mathrm{I}.7.2]{Kol96}$ and $\cite[10.8]{Kol23}$ and the references mentioned in there.

\item
We would like to mention that there is another version (which is in fact easier to comprehend and to use) of the open Chow variety introduced in $\cite[\mathrm{Theorem}\ 21.2]{Har95}$. But since we are working over fields of arbitrary characteristic, we deliberately use the version above in order to be on the safe side.
\end{enumerate}
\end{remark}

\subsection{Some technical preparations}

In this subsection, we do some warm-ups and technical preparations for the proof of Theorem 1.2 in Section 3.

\begin{lemma}
\begin{enumerate}
\item
Let $X,Y,Z$ be varieties and let $f:X\dashrightarrow Y,g:Y\dashrightarrow Z$ be dominant rational maps. Then the composition $g\circ f:X\dashrightarrow Z$ is also dominant.

\item
Let $X,Y$ be two varieties. Let $X_{1},\dots,X_{m};Y_{1},\dots,Y_{n}$ be all of the irreducible components of $X,Y$ respectively. Then $X_{i}\times Y_{j}\subseteq X\times Y\ (1\leq i\leq m,1\leq j\leq n)$ are exactly all of the irreducible components of $X\times Y$.

\item
Let $X,Y$ be two varieties. Let $B\subseteq Y(K)$ be a dense subset of closed points in $Y$, and let $\{A_{y}\subseteq X(K)|\ y\in B\}$ be a collection of sets such that each $A_{y}$ is dense in $X$. Then $\bigsqcup\limits_{y\in B}(A_{y}\times y)\subseteq (X\times Y)(K)$ is dense in $X\times Y$.

\item
Let $X,Y$ be two varieties and let $f:X\dashrightarrow Y,g:X\dashrightarrow Y$ be two rational maps. Suppose $A\subseteq X(K)$ is a dense subset of closed points in $X$ such that $A\subseteq\mathrm{Dom}(f)\cap\mathrm{Dom}(g)$. If $f(x)=g(x)$ for all points $x\in A$, then $f=g$.

\item
Let $X$ be a pure dimensional variety and let $f:X\dashrightarrow X$ be a rational self-map. Let $U\subseteq\mathrm{Dom}(f)$ be an open dense subset of $X$ and let $f_{U}:U\rightarrow X$ be the representative of $f$. Suppose $f_{U}$ is an injective map on $U(K)$, then $f$ is a dominant self-map.
\end{enumerate}
\end{lemma}

\begin{prf}
\begin{enumerate}
\item
This is an immediate consequence of the definition. Let $X^{0},Y^{0},Z^{0}$ be the finite sets of generic points of the irreducible components of $X,Y,Z$ respectively. Then plainly $f(X^{0})=Y^{0}$ and $g(Y^{0})=Z^{0}$ implies that $(g\circ f)(X^{0})=Z^{0}$.

\item
Firstly, every $X_{i}\times Y_{j}$ is an irreducible closed subset of $X\times Y$ since we are working over an algebraically closed field. And we know $\bigcup\limits_{i,j}(X_{i}\times Y_{j})=X\times Y$ by looking ay the closed points. So we only have to prove that pairwise they cannot contain each other. Let $\text{pr}_{1}:X\times Y\rightarrow X$ and $\text{pr}_{2}:X\times Y\rightarrow Y$ be the two projections. Then we finish the proof by noticing $\text{pr}_{1}(X_{i}\times Y_{j})=X_{i}$ and $\text{pr}_{2}(X_{i}\times Y_{j})=Y_{j}$.

\item
Let $F$ be the closure of $\bigsqcup\limits_{y\in B}(A_{y}\times y)$ in $X\times Y$ and let $U=(X\times Y)\backslash F$. Firstly, we have $\bigsqcup\limits_{y\in B}(X\times y)\subseteq F$ because each $A_{y}$ is dense in $X$. So $\text{pr}_{2}(U)\cap B=\emptyset$, where $\text{pr}_{2}:X\times Y\rightarrow Y$ is the projection. But since $\text{pr}_{2}(U)$ is an open subset of $Y$ and $B$ is dense in $Y$, so $\text{pr}_{2}(U)=\emptyset$ and hence $U=\emptyset$. Thus we finish the proof.

\item
Let $U=\mathrm{Dom}(f)\cap\mathrm{Dom}(g)$ be an open dense subset of $X$ and let $f_{U},g_{U}:U\rightarrow Y$ be the representatives of $f,g$ respectively. Then $A\subseteq U(K)$ is dense in $U$ and $f_{U}(x)=g_{U}(x)$ for every $x\in A$. Since $U$ is reduced and $Y$ is separated, we conclude that $f_{U}=g_{U}$ and hence finish the proof. More concretely, the equalizer of $f_{U}$ and $g_{U}$ is a closed subscheme of $U$ containing $A$ because $Y$ is separated. But since $U$ is reduced and $A$ is dense in $U$, this closed subscheme must be $U$ itself. So $f_{U}=g_{U}$.

\item
Let $\eta_{1},\dots,\eta_{m}$ be all of the generic points of $X$. According to our definition of dominant rational maps, we have to prove that $\{f(\eta_{1}),\dots,f(\eta_{m})\}=\{\eta_{1},\dots,\eta_{m}\}$.

Firstly, we prove that $f$ maps generic points to generic points. Assume by contradiction that $f(\eta_1)\notin\{\eta_{1},\dots,\eta_{m}\}$. Let $Y$ be the closure of $f(\eta_1)$ in $X$. Then $Y$ is an irreducible closed subset of $X$ such that $\text{dim}(Y)<\text{dim}(X)$. Moreover, we can find an irreducible open neighborhood $U_{1}$ of $\eta_{1}$ in $U$ such that $f|_{U_1}$ factors through $Y$. Then $f|_{U_1}:U_{1}\rightarrow Y$ is an injective map on $U_{1}(K)$. But since $X$ is pure dimensional, we have $\text{dim}(U_{1})=\text{dim}(X)>\text{dim}(Y)$. So we get a contradiction and hence $f$ maps generic points to generic points.

Secondly, we prove that $f(\eta_i)\neq f(\eta_j)$ for $i\neq j$. Assume by contradiction that $f(\eta_i)=f(\eta_j)=\eta_k$ for some $i\neq j$. Let $X_{t}$ be the irreducible component of $X$ corresponds to $\eta_t$ for each $t=1,\dots,m$; and let $X_{t}^{\circ}=X\backslash(\bigcup\limits_{s\neq t} X_{s})$ be an open dense subset of $X_{t}$ for each $t$. Then $U\cap X_{t}^{\circ}$ is also an open dense subset of $X_{t}$ for each $t$. Now both $f|_{U\cap X_{i}^{\circ}}$ and $f|_{U\cap X_{j}^{\circ}}$ factors through $X_{k}$, and the images $f(U\cap X_{i}^{\circ})$ and $f(U\cap X_{j}^{\circ})$ are both constructible dense subsets of $X_{k}$. By Lemma 2.12 below, we can see that $f(U\cap X_{i}^{\circ})\cap f(U\cap X_{j}^{\circ})$ contains an open dense subset of $X_{k}$. But then $f_{U}$ cannot be an injective map on closed points because $X_{i}^{\circ}\cap X_{j}^{\circ}=\emptyset$. So we get a contradiction and hence prove that $f$ is an injection on the generic points.

We finish the proof by combining the two paragraphs above.
\end{enumerate}
\end{prf}

\begin{lemma}
Let $X,Y,Z$ be varieties and let $f:X\dashrightarrow Y$ be a dominant rational map. Then the rational map $f\times\mathrm{id}_{Z}:X\times Z\dashrightarrow Y\times Z$ is also dominant.
\end{lemma}

\begin{prf}
Let $U=\text{Dom}(f)$ be an open dense subset of $X$ and let $f_{U}:U\rightarrow Y$ be a representative of $f$. Then $U\times Z\subseteq X\times Z$ is an open dense subset (this can be seen by Lemma 2.7(ii)) and $f\times\text{id}_{Z}$ is represented by $f_{U}\times\text{id}_{Z}:U\times Z\rightarrow Y\times Z$. So we may assume that $f$ is a morphism without loss of generality.

Let $X_{i}$ be an irreducible component of $X$. By our definition of dominant, we know $f$ induces a dominant morphism $f_{i}:X_{i}\rightarrow Y_{j(i)}$ where $Y_{j(i)}$ is an irreducible component of $Y$. Now for every irreducible component $Z_{k}$ of $Z$, we can see that $f\times\text{id}_{Z}$ induces the morphism $f_{i}\times\text{id}_{Z_{k}}:X_{i}\times Z_{k}\rightarrow Y_{j(i)}\times Z_{k}$ on the irreducible component $X_{i}\times Z_{k}\subseteq X\times Z$. Using $\cite[\text{Lemma}\ 26.17.5\ \text{(01JT)}]{Stacks}$, we observe that $f_{i}\times\text{id}_{Z_{k}}$ is dominant since $f_{i}$ is dominant. Therefore, we conclude that $f\times\mathrm{id}_{Z}$ is dominant since $f$ is dominant.
\end{prf}

Now we state a technical lemma regarding the domain of definition of a rational map. The proof can be found in $\cite[2.5,\text{Proposition}\ 5]{BLR90}$. Although the original statement contains a smoothness requirement, one can verify that the proof is still valid without that assumption. We also remark that a flat (or more generally, open) morphism satisfies the condition (ii) in our definition of the dominant property, so the composition in the lemma below is valid.

\begin{lemma}
Let $X,X',Y$ be varieties. Let $f:X'\rightarrow X$ be a flat morphism and let $\phi:X\dashrightarrow Y$ be a rational map. Then $\mathrm{Dom}(\phi\circ f)=f^{-1}(\mathrm{Dom}(\phi))$.
\end{lemma}

Next lemma is about the rational group action. Part (i) of the lemma below is $\cite[\text{Lemme}\ 1]{Dem70}$, but we will also include a proof for the convenience of the reader.

\begin{lemma}
Let $G$ be a group variety and let $X$ be an irreducible variety. Let $\rho:G\times X\dashrightarrow X$ be a rational group action of $G$ on $X$.
\begin{enumerate}
\item
For any $g\in G(K)$, we have $\mathrm{Dom}(\rho)\cap(\{g\}\times X)$ is an open dense subset of $\{g\}\times X$. Moreover, the induced rational map $\rho_{g}':X\dashrightarrow X$ given by $X\stackrel{i_{g}}{\hookrightarrow} G\times X\stackrel{\rho}{\dashrightarrow} X$ is dominant, in which $i_{g}:X\hookrightarrow G\times X$ is the closed immersion given by $x\mapsto(g,x)$.

\item
Let $Y,a:G\times Y\rightarrow Y$ and $f:X\dashrightarrow Y$ be a regularization of $\rho$ as in Theorem 2.4(ii). Namely, $Y$ is an irreducible variety, $a$ is a group action of $G$ on $Y$, and $f$ is a $G$-equivariant birational map, i.e. the diagram below is commutative.
\[
\xymatrix{G\times X \ar@{.>}[r]^-{\rho} \ar@{.>}[d]_{\mathrm{id}_{G}\times f} & X \ar@{.>}[d]^{f} \\ G\times Y \ar[r]^-{a} & Y}
\]
Then for every $g\in G(K)$, we have $\rho_{g}'=f^{-1}\circ a_{g}\circ f$ where $a_{g}$ is the automorphism of $Y$ induced by the action $a$ and the point $g$. As a result, $\rho_{g}'$ is birational.
\end{enumerate}
\end{lemma}

\begin{prf}
\begin{enumerate}
\item
Consider the dominant rational map $\phi:G\times X\dashrightarrow X$ given by $\phi=\rho\circ(\sigma\times\mathrm{id}_{X})\circ\tilde{\rho}$, in which $\sigma$ is an automorphism of $G$ given by $h\mapsto gh^{-1}$ and $\tilde{\rho}:G\times X\dashrightarrow G\times X$ is the birational map in Definition 2.2(ii). Then $\phi$ is dominant, and there exists an open dense subset $U\subseteq G\times X$ such that for any $(h,x)\in U(K)$, we have $(h,x)\in\mathrm{Dom}(\tilde{\rho})=\text{Dom}(\rho)$  and $(\sigma\times\mathrm{id}_{X})\circ\tilde{\rho}(h,x)=(gh^{-1},\rho(h,x))\in\mathrm{Dom}(\rho)$.

Now consider the commutative diagram in Definition 2.2(i). For $(h,x)\in U(K)$, we can see that $(gh^{-1},h,x)$ is contained in the domain of definition of the composition map $G\times G\times X\dashrightarrow X$. So by Lemma 2.9, we conclude that $(gh^{-1}\cdot h,x)=(g,x)$ lies in $\text{Dom}(\rho)$. This proves that $\mathrm{Dom}(\rho)\cap(\{g\}\times X)$ is open dense in $\{g\}\times X$. Thus $\rho_{g}'$ is well-defined.

In fact, we have also proved that $\rho_{g}'$ is dominant. Indeed, for each $(h,x)\in U(K)$, we have $(g,x)\in\text{Dom}(\rho)$ and hence $x\in\text{Dom}(\rho_{g}')$. Now $\rho_{g}'(x)=\rho(g,x)=\phi(h,x)$ by the commutative diagram mentioned in the paragraph above. Therefore, we conclude that $\rho_{g}'$ is dominant because $\phi$ is dominant.

\item
Since we have proved that $\rho_{g}'$ is dominant, we may adjust the equality $\rho_{g}'=f^{-1}\circ a_{g}\circ f$ into $f\circ\rho_{g}'=a_{g}\circ f$. We shall just draw a big picture and left the detailed proof to the reader. Let $U\subseteq X$ be $\mathrm{Dom}(f)$ and let $D\subseteq G\times X$ be $\mathrm{Dom}(\rho)$. Notice that $i_{g}^{-1}(D)\cap\rho_{g}'^{-1}(U)\cap U$ is an open dense subset of $X$ thanks to part (i). We also remark that $D\cap\rho^{-1}(U)$ is in fact just $\rho^{-1}(U)$ because $D$ is the domain of definition of $\rho$, but the ``$i_{g}^{-1}(D)$" part in $i_{g}^{-1}(D)\cap\rho_{g}'^{-1}(U)$ should not be dropped.

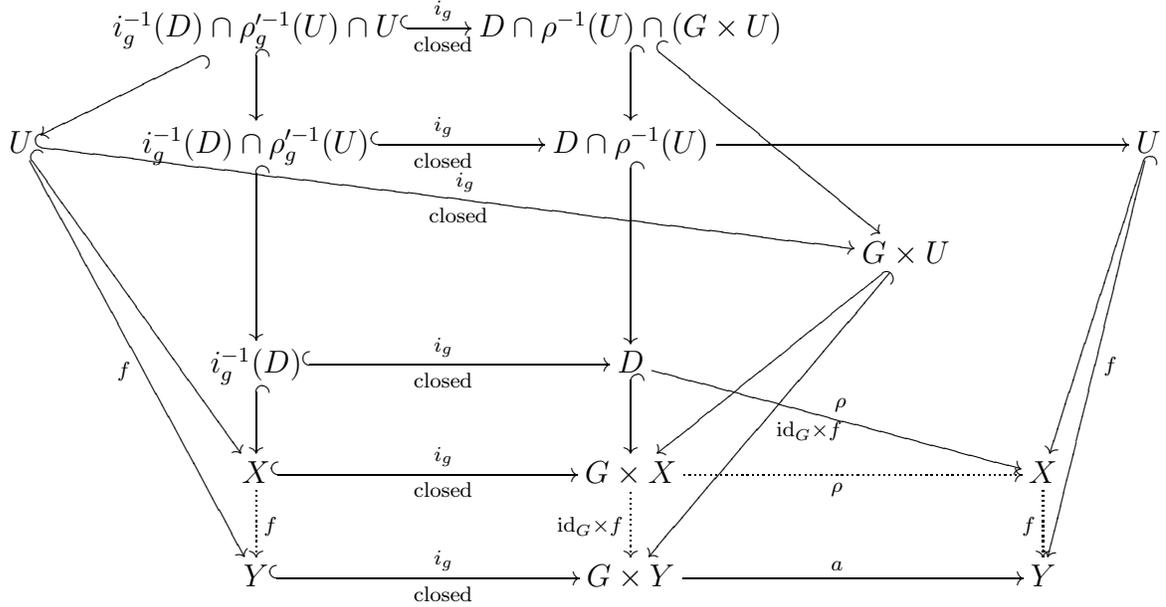
\[
\xymatrix{ & {i_{g}^{-1}(D)\cap\rho_{g}'^{-1}(U)\cap U} \ar@{^{(}->}[r]^{i_{g}}_{\mathrm{closed}} \ar@{^{(}->}[d] \ar@{^{(}->}[ld] & {D\cap\rho^{-1}(U)\cap(G\times U)} \ar@{^{(}->}[d] \ar@{^{(}->}[rdd] & & & \\ U \ar@{^{(}->}[rrrd]^{i_{g}}_{\mathrm{closed}} \ar@{^{(}->}[rddd] \ar[rdddd]_{f} & {i_{g}^{-1}(D)\cap\rho_{g}'^{-1}(U)} \ar@{^{(}->}[r]^{i_{g}}_{\mathrm{closed}} \ar@{^{(}->}[dd] & {D\cap\rho^{-1}(U)} \ar[rrr] \ar@{^{(}->}[dd] & & & U \ar@{^{(}->}[lddd] \ar[ldddd]^{f} \\ & & & {G\times U} \ar@{^{(}->}[ldd] \ar[lddd]^{\mathrm{id}_{G}\times f} & & \\ & {i_{g}^{-1}(D)} \ar@{^{(}->}[r]^{i_{g}}_{\mathrm{closed}} \ar@{^{(}->}[d] & D \ar[rrd]^{\rho} \ar@{^{(}->}[d] & & & \\ & X \ar@{^{(}->}[r]^{i_{g}}_{\mathrm{closed}} \ar@{.>}[d]^{f} & {G\times X} \ar@{.>}[rr]_{\rho} \ar@{.>}[d]_{\mathrm{id}_{G}\times f} & & X \ar@{.>}[d]_{f} & \\ & Y \ar@{^{(}->}[r]^{i_{g}}_{\mathrm{closed}} & G\times Y \ar[rr]^{a} & & Y & }
\]
\end{enumerate}
\end{prf}

We introduce the concept of constructible sets. The definition below is equivalent to the definition of ``locally constructible sets" in $\cite[\text{D\'efinition}\ 9.1.11]{GD61}$ since we are dealing with varieties.

\begin{definition}
Let $X$ be a variety. A \emph{constructible set} in $X$ is a subset of $X$ which is a finite union of locally closed subsets of $X$.
\end{definition}

\begin{lemma}
Let $X$ be a variety and let $A\subseteq X$ be a constructible set. If $A$ is dense in $X$, then $A$ contains an open dense subset of $X$.
\end{lemma}

\begin{prf}
Let $X_{1},\dots,X_{n}$ be the irreducible components of $X$. For each $i=1,\dots,n$, let $X_{i}^{\circ}=X\backslash(\bigcup\limits_{j\neq i} X_{j})$ be an open dense subset of $X_{i}$. Then for each $i$, we know $A\cap X_{i}^{\circ}$ is a constructible dense subset of $X_{i}^{\circ}$. Therefore for each $i$, $A\cap X_{i}^{\circ}$ contains an open dense subset of $X_{i}^{\circ}$ since $X_{i}^{\circ}$ is irreducible. Thus $A$ contains an open dense subset of $X$.
\end{prf}

The set of points in a variety which satisfy a certain property is often constructible. What we need are the following results. The references are $\cite[\text{Th\'er\`eme}\ 9.7.7\text{(i)},\text{Proposition}\ 9.6.1\text{(ii)}]{GD66}$.

\begin{proposition}
Let $X,Y,S$ be varieties. Let $p:X\rightarrow S,q:Y\rightarrow S$ be morphisms, and let $f:X\rightarrow Y$ be an $S$-morphism. Then the following sets are constructible subsets of $S$:
\begin{enumerate}
\item
The subset consists of the points $s\in S$ such that the fiber $X_{s}$ of $p$ is geometrically irreducible.

\item
The subset consists of the points $s\in S$ such that the morphism $f_{s}:X_{s}\rightarrow Y_{s}$ is dominant (in the usual sense, that is, $\overline{f_{s}(X_{s})}=Y_{s}$).
\end{enumerate}
\end{proposition}

\begin{lemma}
Let $f:X\rightarrow S$ be a morphism between varieties and let $U\subseteq X$ be an open dense subset. Then $\{s\in S|\ U_{s}\subseteq X_{s}\text{ is open dense}\}$ is a constructible dense subset of $S$.
\end{lemma}

\begin{prf}
Let $s_{1},\dots,s_{n}$ be all of the generic points of the irreducible components of $S$. Using Proposition 2.13(ii), we see that this subset is constructible in $S$. So we only have to prove that each $s_{i}\in S$ lies in this subset, i.e. $U_{s_{i}}\subseteq X_{s_{i}}$ is an open dense subset. We arbitrarily fix an element $s\in\{s_{1}\dots,s_{n}\}$ and we will prove that $U_{s}\subseteq X_{s}$ is open dense. If $X_{s}=\emptyset$, then there is nothing to prove. So we may assume $X_{s}\neq\emptyset$.

Let $V\subseteq X$ be an open subset. We need to prove that if $V\cap X_{s}\neq\emptyset$, then $V\cap U_{s}=V\cap U\cap X_{s}$ is also nonempty. Equivalently, if there exists an element $v\in V$ such that $f(v)=s$, we need to prove that there exists an element $v'\in V\cap U$ such that $f(v)=s$. Let $V_{0}\subseteq V$ be an irreducible component of $V$ which contains $v$, and let $v_{0}$ be the generic point of $V_{0}$. We claim that $v_{0}\in V\cap U$ and $f(v_{0})=s$, and hence finish the proof.

Firstly, notice that $V\cap U\subseteq V$ is open dense since $U\subseteq X$ is open dense. So $v_{0}$, being the generic point of an irreducible component of $V$, lies in $V\cap U$. Now assume by contradiction that $f(v_{0})\neq s$. Then the closed subset $\overline{f(v_{0})}\subseteq S$ does not contain $s$ because $s$ is the generic point of an irreducible component of $S$. Consider the closed subset $V\cap f^{-1}(\overline{f(v_{0})})$ of $V$. It contains $v_{0}$ but does not contain $v$, which contradicts with our choice of $v_{0}$. So we have proved that $v_{0}\in V\cap U$ and $f(v_{0})=s$. Thus we are done.
\end{prf}

At the end of this subsection, we include a result concerning the graphs of rational maps. We omit the proof since the proof is routine.

\begin{lemma}
Let $X,Y$ be two irreducible varieties.
\begin{enumerate}
\item
Let $f:X\dashrightarrow Y$ be a rational map and let $\Gamma_{f}\subseteq X\times Y$ be the closed graph of $f$. Let $p_{1}:\Gamma_{f}\rightarrow X$, $p_{2}:\Gamma_{f}\rightarrow Y$ be the projections from $\Gamma_{f}$ to $X,Y$ respectively. Then $\Gamma_{f}$ is irreducible, $p_{1}$ is birational and $f=p_{2}\circ p_{1}^{-1}$.

\item
Conversely, let $\Gamma\subseteq X\times Y$ be an irreducible closed subvariety of $X\times Y$ and let $p_{1}:\Gamma\rightarrow X$, $p_{2}:\Gamma\rightarrow Y$ be the projections from $\Gamma$ to $X,Y$ respectively. If $p_{1}$ is birational, then $\Gamma$ is the closed graph of the rational map $p_{2}\circ p_{1}^{-1}:X\dashrightarrow Y$.
\end{enumerate}
\end{lemma}

\begin{corollary}
Let $X,Y$ be two irreducible varieties. Then Lemma 2.15 induces a bijection between the sets $\{f:X\dashrightarrow Y\text{ is a rational map}\}$ and $\{\Gamma\subseteq X\times Y\text{ is an irreducible closed subvariety}$ $\text{such that }p_{1}:\Gamma\rightarrow X\text{ is birational}\}$. Moreover, under this bijection, the birational maps correspond to the irreducible closed subvarieties such that $p_{2}:\Gamma\rightarrow Y$ is also birational.
\end{corollary}

\subsection{A baby case of Theorem 1.2}

In this subsection, we prove Theorem 1.2 in the case that $G$ is finite. There are two reasons that we want to give a short proof of this special case in advance. Firstly, the reader may understand the use of the Weil's regularization theorems through this special case. Secondly, we may then focus on the case that $G$ is infinite in Section 3 and do not let this degenerate case bother us in the main proof.

\begin{proposition}
Let $X$ be an irreducible projective variety and let $G\subseteq\mathrm{Bir}(X)$ be a finite subgroup. Then there exists an irreducible projective variety $Y$ birational to $X$, such that every element of $G$ becomes an automorphism of $Y$ after the birational transformation.
\end{proposition}

\begin{prf}
Since $G$ is a finite group, we may regard $G$ as a finite group variety. Then the birational self-maps in $G$ naturally induce a rational group action $\rho:G\times X\dashrightarrow X$. By Theorem 2.4(ii), there exist an irreducible projective variety $Y$, a group action $a:G\times Y\rightarrow Y$ and a birational map $f:X\dashrightarrow Y$ such that the diagram below is commutative.
\[
\xymatrix{G\times X \ar@{.>}[r]^-{\rho} \ar@{.>}[d]_{\mathrm{id}_{G}\times f} & X \ar@{.>}[d]^{f} \\ G\times Y \ar[r]^-{a} & Y}
\]

Then we finish the proof by Lemma 2.10(ii).
\end{prf}

\section{Proof of Theorem 1.2}

We will prove Theorem 1.2 in this section. Let $X$ be the irreducible projective variety in Theorem 1.2 and fix a very ample line bundle $L\in\text{Pic}(X)$. Let $L_{0}=\text{pr}_{1}^{*}(L)+\text{pr}_{2}^{*}(L)$ be a very ample line bundle on $X\times X$ where $\text{pr}_{1},\text{pr}_{2}:X\times X\rightarrow X$ are the two projections. Let $G\subseteq\text{Bir}(X)$ be the bounded-degree subgroup in Theorem 1.2. Then by definition, $\{(L_{0}^{\text{dim}(X)}\cdot\Gamma_{f})|\ f\in G\}$ is a finite set. We will \textbf{fix the notation} of $X,L$ and $L_{0}$ in this section. We highlight the phrase ``fix the notation" because we will gradually fix many notations in this section, and we want to mark the places where those notations firstly appeared in order to be reader-friendly.

Since we have described our plan of the proof in Introduction, we will just list the steps of the proof here:

\textbf{Step 0:} Construct a quasi-projective variety $S$ such that ``$G\cap S(K)$" (or more precisely, the closed graphs of those elements in $G$) can be naturally viewed as a subset of $S(K)$ which is dense in $S$. Notice that by Corollary 2.16, the closed graphs are one-to-one correspond to the birational maps in $G$.

\textbf{Step 1:} Find an open dense subset $S^{(1)}\subseteq S$, such that every element in $S^{(1)}(K)$ naturally corresponds to a birational self-map of $X$. Moreover, construct a natural ``birational action" of $S$ on $X$, which will later be regularized as a rational group action.

\textbf{Step 2:} Using the group law of $G$, construct a rational group structure on $S$.

\textbf{Step 3:} Using the Weil's regularization theorems, regularize the rational group $S$ and the rational action constructed in Step 1. Then prove that there exists an open dense subset $S^{(2)}\subseteq S^{(1)}$, such that the birational self-maps correspond to elements in $S^{(2)}(K)$ can be simultaneously regularized. That is to say, there exists a birational transformation which turns each of them into an automorphism.

\textbf{Step 4:} Finish the proof by showing that the birational transformation in Step 3 indeed turns every element in $G$ into an automorphism.

~

Now we start with Step 0. We will do the Steps 1--4 in the following subsections 3.1--4.

Using Theorem 2.5 towards $X\times X$ and the very ample line bundle $L_{0}\in\text{Pic}(X\times X)$, we get a quasi-projective variety $C_{\text{dim}(X),d}(X\times X)$ for each positive integer $d$. The closed points in $\bigsqcup\limits_{d=1}^{\infty} C_{\text{dim}(X),d}(X\times X)$ are in one-to-one correspondence with the closed subvarieties of pure dimension $\text{dim}(X)$ in $X\times X$. Since $G$ is of bounded-degree, the subset $\{\Gamma_{f}|\ f\in G\}$ lies in finitely many components of $\bigsqcup\limits_{d=1}^{\infty} C_{\text{dim}(X),d}(X\times X)(K)$. We let $T$ be the closure of $\{\Gamma_{f}|\ f\in G\}$ in $\bigsqcup\limits_{d=1}^{\infty} C_{\text{dim}(X),d}(X\times X)$, which is a quasi-projective variety. We will naturally regard $G$ as a subset of $T(K)$.

Now let $S$ be the union of the irreducible components of maximal dimension in $T$, which is a closed subvariety of $T$. Notice that since $S$ is the closure of an open subset of $T$, we also have that $G\cap S(K)$ is dense in $S$.

The next proposition is an immediate consequence of Theorem 2.5. Notice that $S$ may not be seminormal, but we also do not need it to be seminormal.

\begin{proposition}
Let $S$ be as above. We restrict the universal family on $\bigsqcup\limits_{d=1}^{\infty} C_{\mathrm{dim}(X),d}(X\times X)$ to $S$, and get a closed subvariety $F\subseteq X\times X\times S$. Then the following holds:
\begin{enumerate}
\item
Let $s\in S(K)$ be a closed point. Then the fiber $F_{s}\subseteq X\times X$ has a reduced open dense subset, and $F_{s,\mathrm{red}}$ is tautologically the closed subvariety of $X\times X$ of pure dimension $\mathrm{dim}(X)$ which corresponds to the closed point $s\in\bigsqcup\limits_{d=1}^{\infty} C_{\mathrm{dim}(X),d}(X\times X)(K)$.

\item
If $g\in G\cap S(K)$, then $F_{g,\mathrm{red}}=\Gamma_{g}$ as a closed subvariety of $X\times X$.
\end{enumerate}
\end{proposition}

We will \textbf{fix the notation} of $T,S$ and $F$ in the rest of this section.

\subsection{Step 1}

Recall that in this subsection, our goal is to find an open dense subset $S^{(1)}\subseteq S$ such that every element in $S^{(1)}(K)$ naturally corresponds to a birational self-map of $X$, and construct a natural ``birational action" of $S$ on $X$. The key proposition of this subsection is the result below.

\begin{proposition}
Let $p_{1},p_{2}:F\rightarrow X\times S$ be the two projections as in the diagram below.
\[
\xymatrix{F \ar@/_/[rdd]_{p_{1}} \ar@{^{(}->}[rd] \ar@/^/[rrd]^{p_{2}} & & \\ & X\times X\times S \ar[r]_-{\mathrm{pr}_{2}} \ar[d]^{\mathrm{pr}_{1}} & X\times S \ar[d]^{\mathrm{pr}_{S}} \\ & X\times S \ar[r]_-{\mathrm{pr}_{S}} & S}
\]
Then there exists an open dense subset $S_{1}\subseteq X\times S$ such that both of the induced morphisms $p_{1}^{-1}(S_{1})\rightarrow S_{1}$ and $p_{2}^{-1}(S_{1})\rightarrow S_{1}$ are isomorphisms.
\end{proposition}

We begin with an observation.

\begin{lemma}
\begin{enumerate}
\item
Both of the two morphisms $p_{1},p_{2}:F\rightarrow X\times S$ are surjective.

\item
Let $g\in G\cap S(K)$ be a point and let $p_{1,g},p_{2,g}:F_{g}\rightarrow X$ be the induced morphisms on the fiber. Then there exists an open dense subset $V_{g}\subseteq X$, such that both of the induced morphisms $p_{1,g}^{-1}(V_{g})\rightarrow V_{g}$ and $p_{2,g}^{-1}(V_{g})\rightarrow V_{g}$ are isomorphisms.
\end{enumerate}
\end{lemma}

\begin{prf}
We will just prove the two assertions for $p_{1}$. The same proof can be applied on $p_{2}$.

Let $g\in G\cap S(K)$ be a point, and let $p_{1,g}':F_{g,\text{red}}\rightarrow X$ be the composition of the closed immersion $F_{g,\text{red}}\hookrightarrow F_{g}$ and $p_{1,g}$. Notice that by Proposition 3.1(ii), we know $p_{1,g}'$ is a birational morphism. Moreover, all of the morphisms $p_{1},p_{1,g},p_{1,g}'$ are proper. So in particular, we see that $p_{1,g}'$ is surjective.
\begin{enumerate}
\item
By the discussion above, we can see that $\text{Im}(p_{1})$ is a closed subset of $X\times S$ containing $\bigcup\limits_{g\in G\cap S(K)} (X\times\{g\})$. Since $G\cap S(K)$ is dense in $S$, we conclude that $p_{1}$ is surjective by Lemma 2.7(iii).

\item
Since $F_{g}$ has a reduced open dense subset and $p_{1,g}'$ is birational, we can find an open dense subset $U\subseteq F_{g}$ and an open dense subset $V\subseteq X$ such that $p_{1,g}|_{U}$ factors through an isomorphism $U\stackrel{\sim}\rightarrow V$. Now both $F_{g}$ and $X$ are irreducible, and $\text{dim}(F_{g})=\text{dim}(X)$. So $p_{1,g}(F_{g}\backslash U)$ is a proper closed subset of $X$. One can verify that $V_{g}:=V\backslash(p_{1,g}(F_{g}\backslash U))$ satisfies our requirement.
\end{enumerate}
\end{prf}

The proof of Proposition 3.2 splits into two steps. Firstly, we show that the morphisms $p_{1},p_{2}$ become finite after restricting on some open dense subset.

\begin{lemma}
In the situation of Proposition 3.2, there exists an open dense subset $S_{1}'\subseteq X\times S$ such that both of the induced morphisms $p_{1}^{-1}(S_{1}')\rightarrow S_{1}'$ and $p_{2}^{-1}(S_{1}')\rightarrow S_{1}'$ are finite. Moreover, we can require that $G\cap S(K)\subseteq\mathrm{pr}_{S}(S_{1}')$.
\end{lemma}

\begin{prf}
Let $S_{1}'$ be the subset $\{y\in X\times S|\ \text{dim}(p_{1}^{-1}(y))=\text{dim}(p_{2}^{-1}(y))=0\}$ of $X\times S$. We have seen that $p_{1},p_{2}$ are surjective and proper in Lemma 3.3(i). So $S_{1}'$ is an open subset of $X\times S$. Moreover, the morphisms $p_{1}^{-1}(S_{1}')\rightarrow S_{1}'$ and $p_{2}^{-1}(S_{1}')\rightarrow S_{1}'$ are finite because they are quasi-finite and proper. So it remains to prove that $S_{1}'$ is dense in $X\times S$ and $G\cap S(K)\subseteq\mathrm{pr}_{S}(S_{1}')$.

Using Lemma 3.3(ii), we can see that $\bigcup\limits_{g\in G\cap S(K)}(V_{g}\times\{g\})\subseteq S_{1}'$. Thus we finish the proof by Lemma 2.7(iii).
\end{prf}

The next lemma will be used in the proof of Proposition 3.2. We omit the proof since the proof is routine.

\begin{lemma}
Let $f:X\rightarrow S$ be a finite morphism of schemes and let $S_{0}\subseteq S$ be a closed subscheme. Let $f_{0}:X_{0}:=X\times_{S} S_{0}\rightarrow S_{0}$ be the base change morphism induced by $f$ and let $s_{0}\in S_{0}$ be a point. Then $\mathrm{dim}_{k(s_{0})}((f_{*}\mathcal{O}_{X})_{s_{0}}\otimes_{\mathcal{O}_{S,s_{0}}} k(s_{0}))=\mathrm{dim}_{k(s_{0})}((f_{0,*}\mathcal{O}_{X_{0}})_{s_{0}}\otimes_{\mathcal{O}_{S_{0},s_{0}}} k(s_{0}))$.
\end{lemma}

Now we can prove Proposition 3.2.

\proof[Proof of Proposition 3.2]
The open dense subset $S_{1}\subseteq X\times S$ will be an open subset of $S_{1}'$ in Lemma 3.4. We denote the morphisms $p_{1}^{-1}(S_{1}')\rightarrow S_{1}'$ and $p_{2}^{-1}(S_{1}')\rightarrow S_{1}'$ as $p_{1}'$ and $p_{2}'$ respectively. Then $p_{1}',p_{2}'$ are finite surjective morphisms between varieties. Let $S_{1}$ be the subset $\{y\in S_{1}'|\ \mathrm{dim}_{k(y)}((p_{1,*}'\mathcal{O}_{p_{1}^{-1}(S_{1}')})_{y}\otimes_{\mathcal{O}_{S_{1}',y}} k(y))=\mathrm{dim}_{k(y)}((p_{2,*}'\mathcal{O}_{p_{2}^{-1}(S_{1}')})_{y}\otimes_{\mathcal{O}_{S_{1}',y}} k(y))=1\}$ of $S_{1}'$. Notice that for every point $y\in S_{1}'$, the stalks $(p_{1,*}'\mathcal{O}_{p_{1}^{-1}(S_{1}')})_{y}$ and $(p_{2,*}'\mathcal{O}_{p_{2}^{-1}(S_{1}')})_{y}$ are nonzero. So $S_{1}$ is an open subset of $S_{1}'$.

Denote the morphisms $p_{1}^{-1}(S_{1})\rightarrow S_{1}$ and $p_{2}^{-1}(S_{1})\rightarrow S_{1}$ as $p_{1}''$ and $p_{2}''$ respectively. Then $p_{1}'',p_{2}''$ are finite surjective morphisms between varieties such that $\mathrm{dim}_{k(y)}((p_{1,*}''\mathcal{O}_{p_{1}^{-1}(S_{1})})_{y}\otimes_{\mathcal{O}_{S_{1},y}} k(y))=\mathrm{dim}_{k(y)}((p_{2,*}''\mathcal{O}_{p_{2}^{-1}(S_{1})})_{y}\otimes_{\mathcal{O}_{S_{1},y}} k(y))=1$ for every point $y\in S_{1}$. So one can verify that they are indeed isomorphisms. As a result, it remains to prove that $S_{1}$ is dense in $X\times S$.

For every $g\in G\cap S(K)$, we let $S_{1,g}'$ be the fiber of $S_{1}'$ over $g\in S(K)$. It can be viewed as a nonempty open subset of $X$ because we required that $G\cap S(K)\subseteq\mathrm{pr}_{S}(S_{1}')$. Since $X$ is irreducible, we know that $V_{g}\cap S_{1,g}'$ is an open dense subset of $X$ in which $V_{g}$ is as in Lemma 3.3(ii). Now by Lemma 3.3(ii) and Lemma 3.5, we can see that $\bigcup\limits_{g\in G\cap S(K)}((V_{g}\cap S_{1,g}')\times\{g\})\subseteq S_{1}$. So we finish the proof by Lemma 2.7(iii).
\endproof

We will \textbf{fix the notation} of $S_{1}$ in the rest of this section. Notice that $G\cap S(K)\subseteq\mathrm{pr}_{S}(S_{1})$ also holds. 

~

Now we can describe the ``birational action" of $S$ on $X$. Let $\sigma_{1},\sigma_{2}:X\times S\dashrightarrow X\times S$ be the rational maps determined by the composite morphisms $S_{1}\stackrel{\sim}\rightarrow p_{1}^{-1}(S_{1})\hookrightarrow F\stackrel{p_{2}}\rightarrow X\times S$ and $S_{1}\stackrel{\sim}\rightarrow p_{2}^{-1}(S_{1})\hookrightarrow F\stackrel{p_{1}}\rightarrow X\times S$ respectively. Then $\sigma_{1},\sigma_{2}$ are $S$-rational maps. We let $\rho:=\text{pr}_{X}\circ\sigma_{1}:X\times S\dashrightarrow X$ be our ``birational action" (notice that rational maps can always composite with morphisms in this way). We will later prove that (see Proposition 3.8) both $\sigma_{1}$ and $\sigma_{2}$ are dominant and in fact they are inverse to each other, i.e. $\sigma_{1}\circ\sigma_{2}=\sigma_{2}\circ\sigma_{1}=\text{id}_{X\times S}$. But now, we start by analyzing the behavior of $\sigma_{1}$ and $\sigma_{2}$ on the fibers.

\begin{definition}
Let $s\in S(K)$ be a closed point contained in $\mathrm{pr}_{S}(S_{1})$. Then $S_{1,s}$ is an open dense subset of $X$. There are two equivalent ways to define the rational maps $\sigma_{1,s},\sigma_{2,s}:X\dashrightarrow X$. We will focus on $\sigma_{1,s}$ below and the situation for $\sigma_{2,s}$ is just the same:
\begin{enumerate}
\item
Define $\sigma_{1,s}$ as the base change of $\sigma_{1}$ directly, i.e. let $\sigma_{1,s}$ be the self-map of $X$ induced by the composition $S_{1,s}\stackrel{\sim}\rightarrow p_{1}^{-1}(S_{1})_{s}\hookrightarrow F_{s}\stackrel{p_{2,s}}\longrightarrow X$.

\item
Let $\sigma_{1,s}$ be induced by the action $\rho$ in the way as in Lemma 2.10(i). More specifically, let $\sigma_{1,s}$ be the self-map of $X$ induced by the composition $S_{1,s}\hookrightarrow S_{1}\stackrel{\rho}\rightarrow X$ where $S_{1,s}\hookrightarrow S_{1}$ is the natural closed immersion.
\end{enumerate}
One can verify that the two ways above indeed give the same $\sigma_{1,s}$ (resp. $\sigma_{2,s}$).
\end{definition}

Notice that $\mathrm{pr}_{S}(S_{1})$ is an open dense subset of $S$. So we have defined the notions $\sigma_{1,s},\sigma_{2,s}$ on ``general closed points" of $S$. Moreover, one can verify that for $g\in G\cap S(K)\subseteq\mathrm{pr}_{S}(S_{1})$, we have $\sigma_{1,g}=\sigma_{2,g}^{-1}=g$ as a self-map of $X$. 

We also denote $\rho_{s}:=\sigma_{1,s}$ for $s\in\text{pr}_{S}(S_{1})(K)$. So $\rho_{s}$ and $\sigma_{1,s}$ are just two names of one object (but $\sigma_{2,s}$ only has one name). We will use different names in different places to indicate which one of the two equivalent definitions is in our mind at that place. We will \textbf{fix the notation} of $\sigma_{1},\sigma_{2},\rho$ and $\sigma_{1,s}=\rho_{s},\sigma_{2,s}$ as above in the rest of this section.

~

Now we can describe the open dense subset $S^{(1)}\subseteq S$ which appears at the beginning of this subsection. We apply Proposition 2.13(i) to the structure morphism $F\rightarrow S$ (see the diagram in Proposition 3.2) and get that $\{s\in S|\ F_{s}\text{ is geometrically irreducible}\}$ is a constructible subset of $S$. Using Proposition 3.1(ii), we see that $G\cap S(K)$ is contained in this set. So this set is a constructible dense subset of $S$, and hence contains an open dense subset of $S$ by Lemma 2.12. As a result, we can choose an open dense subset $S^{(1)}$ of $S$ such that $S^{(1)}\subseteq\{s\in S|\ F_{s}\text{ is geometrically irreducible}\}\cap\mathrm{pr}_{S}(S_{1})$. We \textbf{fix the notation} of $S^{(1)}$ in the rest of this section. We remark that unfortunately, we do not know whether we can guarantee $G\cap S(K)\subseteq S^{(1)}(K)$.

The next proposition shows that the closed points in $S^{(1)}$ naturally correspond to some birational self-maps of $X$.

\begin{proposition}
Let $s\in S^{(1)}(K)$ be a closed point. Then $F_{s,\mathrm{red}}\subseteq X\times X$ is an irreducible closed subvariety. Let $p_{1,s}',p_{2,s}':F_{s,\mathrm{red}}\rightarrow X$ be the two projections constructed as in the proof of Lemma 3.3. Then we have:
\begin{enumerate}
\item
Both $p_{1,s}'$ and $p_{2,s}'$ are birational morphisms.

\item
We have $\sigma_{1,s}=p_{2,s}'\circ p_{1,s}'^{-1}$ and $\sigma_{2,s}=p_{1,s}'\circ p_{2,s}'^{-1}$. Therefore, $\sigma_{1,s}=\rho_{s}$ is a birational self-map of $X$ and $\sigma_{2,s}$ is its inverse.

\item
The irreducible closed subvariety $F_{s,\mathrm{red}}\subseteq X\times X$ is the closed graph of $\sigma_{1,s}=\rho_{s}$.

\item
The closed point $s\in S(K)\subseteq\bigsqcup\limits_{d=1}^{\infty} C_{\mathrm{dim}(X),d}(X\times X)(K)$ corresponds to the closed graph of $\rho_{s}$.
\end{enumerate}
\end{proposition}

\begin{prf}
\begin{enumerate}
\item
We will prove this for $p_{1,s}'$. The proof for $p_{2,s}'$ is just the same.

Firstly, notice that $p_{1,s}:F_{s}\rightarrow X$ naturally induces an isomorphism $p_{1}^{-1}(S_{1})_{s}\stackrel{\sim}\rightarrow S_{1,s}$. Since $s\in\mathrm{pr}_{S}(S_{1})(K)$, we know $S_{1,s}$ is an nonempty open subscheme of $X$. As a result, we see that $p_{1}^{-1}(S_{1})_{s}$ is an integral open subscheme of $F_{s}$. Hence the open immersion $p_{1}^{-1}(S_{1})_{s}\hookrightarrow F_{s}$ factors through $F_{s,\text{red}}$, and $p_{1}^{-1}(S_{1})_{s}$ can be viewed as an open dense subset of $F_{s,\text{red}}$ because $F_{s,\text{red}}$ is irreducible. Now $p_{1,s}'|_{p_{1}^{-1}(S_{1})_{s}}$ is an open immersion onto $S_{1,s}\subseteq X$, so $p_{1,s}'$ is a birational morphism.

\item
This follows from the definition of $\sigma_{1,s},\sigma_{2,s}$ directly. See also the paragraph above.

\item
This follows from part (ii) and Lemma 2.15(ii) directly.

\item
This tautologically follows from part (iii). See Proposition 3.1(i).
\end{enumerate}
\end{prf}

Now we have finished half of the Step 1. We will prove the promised properties of $\sigma_{1},\sigma_{2}$ and thus finish the Step 1.

\begin{proposition}
\begin{enumerate}
\item
Both $\sigma_{1}$ and $\sigma_{2}$ are dominant self-maps of $X\times S$.

\item
We have $\sigma_{1}\circ\sigma_{2}=\sigma_{2}\circ\sigma_{1}=\mathrm{id}_{X\times S}$. So $\sigma_{1},\sigma_{2}$ are birational self-maps of $X\times S$.
\end{enumerate}
\end{proposition}

\begin{prf}
\begin{enumerate}
\item
We will prove this assertion for $\sigma_{1}$. The proof for $\sigma_{2}$ is just the same.

By definition, we have $S_{1}\subseteq\text{Dom}(\sigma_{1})$. Since $\sigma_{1}$ is a self-map of $X\times S$, we only have to prove that $\sigma_{1}(S_{1})$ is dense in $X\times S$ in order to prove $\sigma_{1}$ is dominant (recall our definition about dominant rational maps in Section 1). Let $s\in S^{(1)}(K)$ be a closed point. Then by Proposition 3.7(ii), there exists an open dense subset $V_{s}$ of $X$ such that $V_{s}\subseteq\sigma_{1,s}(S_{1,s})$. Then one can see that $\bigcup\limits_{s\in S^{(1)}(K)}(V_{s}(K)\times\{s\})\subseteq\sigma_{1}(S_{1})$. So we finish the proof by Lemma 2.7(iii).

\item
We will prove $\sigma_{2}\circ\sigma_{1}=\mathrm{id}_{X\times S}$. The proof of $\sigma_{1}\circ\sigma_{2}=\mathrm{id}_{X\times S}$ is just the same.

By Lemma 2.7(iv), we only need to find a subset $A\subseteq S_{1}(K)$, such that:
\begin{enumerate}
\item
$A$ is dense in $X\times S$.

\item
For every element $a\in A$, we have $\sigma_{1}(a)\in S_{1}(K)$ and $\sigma_{2}(\sigma_{1}(a))=a$.
\end{enumerate}
Now let $s\in S^{(1)}(K)$ be a closed point. Still by Proposition 3.7(ii), we can see that there exists an open dense subset $U_{s}$ of $X$ contained in $S_{1,s}$ such that for every element $x\in U_{s}(K)$, we have $\sigma_{1,s}(x)\in S_{1,s}(K)$ and $\sigma_{2,s}(\sigma_{1,s}(x))=x$. Then one can verify that $A:=\bigcup\limits_{s\in S^{(1)}(K)}(U_{s}(K)\times\{s\})$ satisfies our requirements and hence we are done.
\end{enumerate}
\end{prf}

\subsection{Step 2}

In this subsection, our goal is to construct a rational group structure on $S$ which is compatible with the original group law of $G$. The idea is easy: we will find an appropriate closed subvariety $\Gamma\subseteq X\times X\times S\times S$, which can induce the right rational map $m:S\times S\dashrightarrow S$ via the universal property of the open Chow variety. The closed subvariety $\Gamma$ will be defined as a closed graph of an $(S\times S)$-rational map $X\times S\times S\dashrightarrow X\times S\times S$.

In order to make things clearly, we distinguish the two $S$ in $X\times S\times S$ and write it as $X\times S_{(1)}\times S_{(2)}$. This is a temporary notation and will only appears in this paragraph. Let $f_{1}$ be the birational map $\sigma_{1}\times\text{id}_{S_{(2)}}:X\times S_{(1)}\times S_{(2)}\dashrightarrow X\times S_{(1)}\times S_{(2)}$ in which $\sigma_{1}$ is regarded as a birational self-map of $X\times S_{(1)}$; and let $f_{2}$ be the birational map $\sigma_{1}\times\text{id}_{S_{(1)}}:X\times S_{(1)}\times S_{(2)}\dashrightarrow X\times S_{(1)}\times S_{(2)}$ in which $\sigma_{1}$ is regarded as a birational self-map of $X\times S_{(2)}$. Then both $f_{1}$ and $f_{2}$ are birational $(S\times S)$-self-maps of $X\times S\times S$.

The next lemma directly follows from the definition of $f_{1}$ and $f_{2}$. We omit its proof since the proof is routine.

\begin{lemma}
Let $(s_{1},s_{2})\in(S^{(1)}\times S^{(1)})(K)$ be a closed point.
\begin{enumerate}
\item
The fiber $\mathrm{Dom}(f_{1})_{(s_{1},s_{2})}$ is an open dense subset of $X$. Hence we can define $f_{1,(s_{1},s_{2})}:X\dashrightarrow X$ via the morphism $\mathrm{Dom}(f_{1})_{(s_{1},s_{2})}\rightarrow X$. Then $f_{1,(s_{1},s_{2})}=\sigma_{1,s_{1}}=\rho_{s_{1}}$.

\item
Similar result holds for $f_{2}$, while the last equation turns to $f_{2,(s_{1},s_{2})}=\sigma_{1,s_{2}}=\rho_{s_{2}}$.
\end{enumerate}
\end{lemma}

Let $f:=f_{1}\circ f_{2}$ which is also a birational $(S\times S)$-self-map of $X\times S\times S$. Then a similar result also holds for $f$, which can be easily deduced using the fact that $\rho_{s_{1}}$ and $\rho_{s_{2}}$ are dominant. So we also omit the proof.

\begin{lemma}
Let $(s_{1},s_{2})\in(S^{(1)}\times S^{(1)})(K)$ be a closed point. Then the fiber $\mathrm{Dom}(f)_{(s_{1},s_{2})}$ is an open dense subset of $X$, and the induced rational map $f_{(s_{1},s_{2})}:X\dashrightarrow X$ equals to $\rho_{s_{1}}\circ\rho_{s_{2}}$.
\end{lemma}

Now we define $\Gamma\subseteq X\times X\times S\times S$ as the closed graph of the $(S\times S)$-rational map $f$. See the diagram below. We denote $U$ as the open dense subset $\text{Dom}(f)\subseteq X\times S\times S$, and the underlines below $X$ are used to indicate that those morphisms are obtained by timing with \emph{that} $\text{id}_{X}$. The closed subvariety $\Gamma$ is the closure of the locally closed subvariety $U\subseteq X\times X\times S\times S$. Notice that the open immersion $U\hookrightarrow\Gamma$ at the top of this diagram makes $U$ into an open dense subset of $\Gamma$.
\[
\xymatrix{U \ar@{^{(}->}[r] \ar@/^3pc/[rrrd]^{f} \ar@{^{(}->}[rd]^{\text{closed}} \ar@/_/[rdd]_{\text{id}} & \Gamma \ar@{^{(}->}[rd]^{\text{closed}} & & \\ & U\times\underline{X} \ar@{^{(}->}[r] \ar[d]^{\mathrm{pr}_{U}} & X\times\underline{X}\times S\times S \ar[r]^-{\mathrm{pr}_{2}} \ar[d]^{\mathrm{pr}_{1}} & \underline{X}\times S\times S \ar[d]^{\mathrm{pr}_{S\times S}} \\ & U \ar@{^{(}->}[r] & X\times S\times S \ar[r]^-{\mathrm{pr}_{S\times S}} & S\times S}
\]

We denote $p:\Gamma\rightarrow S\times S$ as the surjective structure morphism obtained by the diagram above. Using the generic flatness theorem $\cite[\text{Proposition}\ 29.27.2\ \text{(052B)}]{Stacks}$, Lemma 2.12 and Lemma 2.14, we can see that there exists an open dense subset $Q\subseteq S\times S$ such that:
\begin{enumerate}
\item
$Q$ is a nonsingular open subvariety of $S^{(1)}\times S^{(1)}$.

\item
The morphism $p^{-1}(Q)\rightarrow Q$ induced by $p$ is flat.

\item
For each point $q\in Q$, the fiber $U_{q}\subseteq X_{k(q)}$ can be viewed as an open dense subset of $\Gamma_{q}$.
\end{enumerate}

The next lemma is an immediate consequence of Lemma 3.10 and the part (iii) above.

\begin{lemma} 
\begin{enumerate}
\item
For each closed point $q=(s_{1},s_{2})\in Q(K)$, the closed subvariety $\Gamma_{q,\mathrm{red}}\subseteq X\times X$ is the closed graph of the birational map $\rho_{s_{1}}\circ\rho_{s_{2}}:X\dashrightarrow X$.

\item
For each point $q\in Q$, the fiber $\Gamma_{q}\subseteq(X\times X)_{k(q)}$ is an irreducible closed subscheme of dimension equals to $\mathrm{dim}(X)$ and contains a geometrically reduced open dense subset (which is $U_{q}$).
\end{enumerate}
\end{lemma}

Now we can introduce the rational group law on $S$. It is constructed via the functoriality of the open Chow variety.

\begin{lemma}
There exists a morphism $m':Q\rightarrow\bigsqcup\limits_{d=1}^{\infty} C_{\mathrm{dim}(X),d}(X\times X)$ such that for each point $q=(s_{1},s_{2})\in Q(K)$, the closed point $m'(q)\in\bigsqcup\limits_{d=1}^{\infty} C_{\mathrm{dim}(X),d}(X\times X)(K)$ corresponds to the closed graph of $\rho_{s_{1}}\circ\rho_{s_{2}}:X\dashrightarrow X$.
\end{lemma}

\begin{prf}
Consider the closed subvariety $p^{-1}(Q)\subseteq X\times X\times Q$. By Lemma 3.11(ii), this flat family over $Q$ satisfies every condition about the family in Theorem 2.5, except one thing that the degree of the fibers may not be constant. But since this family is flat, we know that the degree of the fibers (calculated by intersecting with $L_{0}\in\text{Pic}(X\times X)$) is a locally constant function on $Q$. So we can write $Q=\bigsqcup\limits_{d=1}^{\infty} Q_{d}$ where $Q_{d}$ is the open subset of $Q$ containing those points at which the fiber has degree $d$. Hence by the functoriality, we can get a desired morphism $Q_{d}\rightarrow C_{\mathrm{dim}(X),d}(X\times X)$ for each $d$ and then assemble them into a morphism $m':Q\rightarrow\bigsqcup\limits_{d=1}^{\infty} C_{\mathrm{dim}(X),d}(X\times X)$. Notice that each variety $Q_{d}$ is nonsingular and thus seminormal.

The assertion regarding the closed points $q\in Q(K)$ is guaranteed by Lemma 3.11(i).
\end{prf}

\begin{lemma}
The morphism $m':Q\rightarrow\bigsqcup\limits_{d=1}^{\infty} C_{\mathrm{dim}(X),d}(X\times X)$ in the lemma above factors through $T\subseteq\bigsqcup\limits_{d=1}^{\infty} C_{\mathrm{dim}(X),d}(X\times X)$.
\end{lemma}

\begin{prf}
Since $Q$ is reduced, we only need to prove that $m'(Q)\subseteq T$. Moreover, we only need to find a dense subset $A\subseteq Q$ such that $m'(A)\subseteq T$ since $T$ is closed in $\bigsqcup\limits_{d=1}^{\infty} C_{\mathrm{dim}(X),d}(X\times X)$.

Recall that $G\cap S(K)$ is dense in $S$. So $(G\cap S(K))\times(G\cap S(K))$ is a dense subset of closed points in $S\times S$ by Lemma 2.7(iii). As a result, the set $A:=((G\cap S(K))\times(G\cap S(K)))\cap Q(K)$ is a dense subset of $Q$. For each point $(g_{1},g_{2})\in A$, we have $\rho_{g_{1}}=g_{1}$ and $\rho_{g_{2}}=g_{2}$ as birational self-maps of $X$ by Proposition 3.1(ii) and Proposition 3.7(iii) (and the one-to-one correspondence in Corollary 2.16). So $m'((g_{1},g_{2}))$ corresponds to the closed graph of $\rho_{g_{1}}\circ\rho_{g_{2}}=g_{1}g_{2}$ and thus we can say that $m'((g_{1},g_{2}))=g_{1}g_{2}\in G\subseteq T(K)$. Therefore, we find a dense subset $A$ of $Q$ such that $m'(A)\subseteq T$ and hence finish the proof.
\end{prf}

Abusing notation, we may regard $m'$ as a morphism $Q\rightarrow T$ and also denote this morphism as $m'$.

\begin{lemma}
The morphism $m':Q\rightarrow T$ factors through $S$.
\end{lemma}

\begin{prf}
As in the proof of the lemma above, we only have to show that $m'(Q)\subseteq S$. Assume by contradiction that $m'(Q)\cap(T\backslash S)\neq\emptyset$. Since $T\backslash S$ is an open subset of $T$, there exists a generic point $\eta$ of $Q$ such that $m'(\eta)\in T\backslash S$. Let $V$ be the closure of $m'(\eta)$ in $T$. Then by the definition of $S$, we know that $V$ is an irreducible closed subset of $T$ such that $\text{dim}(V)<\text{dim}(T)$. Moreover, we can find an irreducible open neighborhood $U$ of $\eta$ in $Q$ such that $m'|_{U}$ factors through $V$.

Since $Q$ is an open dense subset of $S\times S$ and $S$ is pure dimension of $\text{dim}(T)$, we can see that $S\times S$ and hence $Q$ is pure dimensional of $2\cdot\text{dim}(T)$ by Lemma 2.7(ii). In particular, we have $\text{dim}(U)=2\cdot\text{dim}(T)$. Now consider the morphism $\Phi:U\rightarrow V\times S$ given by $m'|_{U}:U\rightarrow V$ and $U\hookrightarrow Q\hookrightarrow S\times S\stackrel{\text{pr}_2}\rightarrow S$. Then for each point $u=(s_{1},s_{2})\in U(K)$, we have $\Phi(u)=(m'(u),s_{2})$ where $m'(u)\in V\subseteq T$ is the closed point of $\bigsqcup\limits_{d=1}^{\infty} C_{\mathrm{dim}(X),d}(X\times X)$ which corresponds to the closed graph of $\rho_{s_{1}}\circ\rho_{s_{2}}$. As a result, we know that $\Phi$ is an injective map on $U(K)$.

However, the dimension $\text{dim}(U)=2\cdot\text{dim}(T)>\text{dim}(V)+\text{dim}(S)$. So we get a contradiction and thus conclude that $m'$ factors through $S$.
\end{prf}

Now we conclude what we have gained through the discussion above in part (i) of the proposition below. Replacing the birational $(S\times S)$-self map $f=f_{1}\circ f_{2}$ of $X\times S\times S$ by $f_{1}\circ f_{2}^{-1}$ and $f_{1}^{-1}\circ f_{2}$ respectively and applying just the same argument as above, one then gets the conclusion of part (ii) and (iii). 

\begin{proposition}
\begin{enumerate}
\item
There exists an open dense subset $Q\subseteq S^{(1)}\times S^{(1)}$ and a morphism $m:Q\rightarrow S$ such that for each point $q=(s_{1},s_{2})\in Q(K)$, we have $m(q)$ is the closed point in $\bigsqcup\limits_{d=1}^{\infty} C_{\mathrm{dim}(X),d}(X\times X)$ which corresponds to the closed graph of $\rho_{s_{1}}\circ\rho_{s_{2}}:X\dashrightarrow X$.

\item
There exists an open dense subset $Q_{1}\subseteq S^{(1)}\times S^{(1)}$ and a morphism $m_{1}:Q_{1}\rightarrow S$ such that for each point $q=(s_{1},s_{2})\in Q_{1}(K)$, we have $m_{1}(q)$ is the closed point in $\bigsqcup\limits_{d=1}^{\infty} C_{\mathrm{dim}(X),d}(X\times X)$ which corresponds to the closed graph of $\rho_{s_{1}}\circ\rho_{s_{2}}^{-1}:X\dashrightarrow X$.

\item
There exists an open dense subset $Q_{2}\subseteq S^{(1)}\times S^{(1)}$ and a morphism $m:Q_{2}\rightarrow S$ such that for each point $q=(s_{1},s_{2})\in Q_{2}(K)$, we have $m_{2}(q)$ is the closed point in $\bigsqcup\limits_{d=1}^{\infty} C_{\mathrm{dim}(X),d}(X\times X)$ which corresponds to the closed graph of $\rho_{s_{1}}^{-1}\circ\rho_{s_{2}}:X\dashrightarrow X$.
\end{enumerate}
Notice that $Q,Q_{1},Q_{2}$ are also open dense subsets of $S\times S$. We emphasize that they are contained in $S^{(1)}\times S^{(1)}$ because we need to guarantee that the notions $\rho_{s_{1}},\rho_{s_{2}}$ are well-defined.
\end{proposition}

We regard the morphisms $m,m_{1},m_{2}$ in the proposition above as rational maps $m,m_{1},m_{2}:S\times S\dashrightarrow S$. We will show that $m$ is a rational group law that satisfies the conditions in Definition 2.1. The maps $m_{1},m_{2}$ will play an auxiliary role in the proof that $m$ satisfies Definition 2.1(ii). We will write the maps $m,m_{1},m_{2}$ as $S\times S\dashrightarrow S:(v,w)\mapsto vw,vw^{-1},v^{-1}w$ respectively.

We \textbf{fix the notation} of $Q,Q_{1},Q_{2}$ and $m,m_{1},m_{2}$ as above in the rest of this section.

~

Now we prove that $m$ is a rational group law. We start with an easy observation. It is an immediate consequence of Proposition 3.1(i), Proposition 3.7(iii) and the one-to-one correspondence in Corollary 2.16. So we omit the proof.

\begin{lemma}
Let $\Phi:S\times S\dashrightarrow S\times S$ be the rational map given by $(v,w)\mapsto(vw,w)$. Then $Q\subseteq\mathrm{Dom}(\Phi)$ and $\Phi$ is an injective map on $Q(K)$. Similar result holds for:
\begin{enumerate}
\item
$\Phi':S\times S\dashrightarrow S\times S,(v,w)\mapsto(vw^{-1},w)$ and $Q_{1}$.

\item
$\Psi:S\times S\dashrightarrow S\times S,(v,w)\mapsto(v,vw)$ and $Q$.

\item
$\Psi':S\times S\dashrightarrow S\times S,(v,w)\mapsto(v,v^{-1}w)$ and $Q_{2}$.
\end{enumerate}
\end{lemma}

We also \textbf{fix the notation} of $\Phi,\Phi',\Psi,\Psi'$ as above in the rest of this section. The next lemma is an immediate consequence of Lemma 2.7(ii)(v) since $S$ and hence $S\times S$ is of pure dimensional.

\begin{lemma}
All of the rational maps $\Phi,\Phi',\Psi,\Psi'$ are dominant.
\end{lemma}

Now we can prove that $m$ satisfies Definition 2.1(ii).

\begin{lemma}
We have $\Phi\circ\Phi'=\Phi'\circ\Phi=\Psi\circ\Psi'=\Psi'\circ\Psi=\mathrm{id}_{S\times S}$. So $\Phi,\Phi',\Psi,\Psi'$ are birational self-maps of $S\times S$ and hence $m$ satisfies the condition in Definition 2.1(ii).
\end{lemma}

\begin{prf}
We will prove that $\Phi'\circ\Phi=\mathrm{id}_{S\times S}$. The proof for the other three equalities are same.

We regard $\Phi,\Phi'$ as morphisms $\Phi|_{Q}:Q\rightarrow S\times S$ and $\Phi'|_{Q_1}:Q_{1}\rightarrow S\times S$. By the lemma above, we know that $(\Phi|_{Q})^{-1}(Q_1)$ is an open dense subset of $S\times S$. We prove that for every closed point $q=(s_1,s_2)\in(\Phi|_{Q})^{-1}(Q_1)(K)$, we have $\Phi'(\Phi(q))=q$. If so, then we are done by Lemma 2.7(iv).

Firstly, we have $\Phi(q)=(m(s_1,s_2),s_2)\in Q_1(K)$. We know that $m(s_1,s_2)$ is the closed point in the open Chow variety which corresponds to the closed graph of $\rho_{s_1}\circ\rho_{s_2}$. But by Proposition 3.7(iii), we have that $m(s_1,s_2)$ also corresponds to the closed graph of $\rho_{m(s_1,s_2)}$. So by the one-to-one correspondence in Corollary 2.16, we know $\rho_{m(s_1,s_2)}=\rho_{s_1}\circ\rho_{s_2}$.

Now $\Phi'(\Phi(q))=(m_{1}(m(s_1,s_2)),s_2)$ in which $m_{1}(m(s_1,s_2))$ is the closed point in the open Chow variety which corresponds to the closed graph of $\rho_{m(s_1,s_2)}\circ\rho_{s_2}^{-1}=\rho_{s_1}$. So we have $m_{1}(m(s_1,s_2))=s_1$ and hence $\Phi'(\Phi(q))=q$ for every closed point $q=(s_1,s_2)\in(\Phi|_{Q})^{-1}(Q_1)(K)$. Thus we finish the proof.
\end{prf}

Now we prove that $m$ satisfies Definition 2.1(i). Firstly, notice that $m$ is a dominant rational map because $m=\text{pr}_{1}\circ\Phi$ where $\text{pr}_{1}$ is the projection from $S\times S$ to $S$. So taking Lemma 2.8 into account, we know that all of the arrows in the diagram of Definition 2.1(i) are dominant rational maps. As a result, the compositions are well-defined.

\begin{lemma}
The diagram
\[
\xymatrix{S\times S\times S \ar@{.>}[r]^-{\mathrm{id}_{S}\times m} \ar@{.>}[d]_{m\times\mathrm{id}_{S}} & S\times S \ar@{.>}[d]^{m} \\ S\times S \ar@{.>}[r]^{m} & S}
\]
is commutative and hence $m$ satisfies the condition in Definition 2.1(i).
\end{lemma}

\begin{prf}
Let $p_{1}:Q\times S\stackrel{m\times\mathrm{id}_{S}}\longrightarrow S\times S$ and $p_{2}:S\times Q\stackrel{\mathrm{id}_{S}\times m}\longrightarrow S\times S$ be the representatives of the rational maps $m\times\mathrm{id}_{S}$ and $\mathrm{id}_{S}\times m$ respectively. As we have seen in the paragraph above, the maps $p_1,p_2$ and $m$ are dominant. So $p_{1}^{-1}(Q)\cap p_{2}^{-1}(Q)$ is an open dense subset of $S\times S\times S$.

Let $q=(s_1,s_2,s_3)$ be a closed point in $(p_{1}^{-1}(Q)\cap p_{2}^{-1}(Q))(K)$. Using the argument as in the proof of Lemma 3.18, we can see that both of the points $m(p_1(q)),m(p_2(q))\in S(K)$ are the closed point in the open Chow variety which corresponds to the closed graph of $\rho_{s_1}\circ\rho_{s_2}\circ\rho_{s_3}$. So we conclude that $m(p_1(q))=m(p_2(q))$ for every point $q\in(p_{1}^{-1}(Q)\cap p_{2}^{-1}(Q))(K)$. Therefore we finish the proof by Lemma 2.7(iv).
\end{prf}

According Lemma 3.18 and Lemma 3.19, we can make a conclusion in the proposition below.

\begin{proposition}
The rational map $m:S\times S\dashrightarrow S$ makes $S$ into a rational group in the sense of Definition 2.1.
\end{proposition}

There is also one thing that we should do in this subsection. As we have described the ``group law" on $S$, we will prove that this group law is compatible with the ``group action" $\rho:X\times S\dashrightarrow X$.

\begin{lemma}
The diagram
\[
\xymatrix{S\times S\times X \ar@{.>}[r]^-{\mathrm{id}_{S}\times\rho} \ar@{.>}[d]_{m\times\mathrm{id}_{X}} & S\times X \ar@{.>}[d]^{\rho} \\ S\times X \ar@{.>}[r]^{\rho} & X}
\]
is commutative. Notice that we interchange the position of $X$ and $S$ in $X\times S$ and regard $\rho$ as a rational map $S\times X\dashrightarrow X$ because we want a left action as a matter of convention.
\end{lemma}

\begin{prf}
Firstly, one can verify that all of the arrows in this diagram are dominant rational maps as usual. Recall that for every point $s\in S^{(1)}(K)$, we have an open dense subset $S_{1,s}$ of $X$. For every point $x\in S_{1,s}(K)$, we have that $x\in\text{Dom}(\rho_s),(s,x)\in\text{Dom}(\rho)$ and $\rho_{s}(x)=\rho((s,x))$.

We regard $m$ as a morphism $Q\rightarrow S$ and denote $R$ as $m^{-1}(S^{(1)})$. Since $m$ is dominant, we know that $R\subseteq Q$ is an open dense subset of $S\times S$. Now for each point $r=(s_1,s_2)\in R(K)$, consider the open dense subset $A_{(s_1,s_2)}:=\rho_{s_2}^{-1}(S_{1,s_1})\cap S_{1,m(s_1,s_2)}$ of $X$. Here we regard $\rho_{s_2}$ as a morphism $S_{1,s_2}\rightarrow X$. Notice that since the points $s_1,s_2,m(s_1,s_2)$ all lie in $S^{(1)}(K)$, this notion is well-defined. Moreover, $A_{(s_1,s_2)}$ is indeed an open dense subset of $X$ because the self-maps $\rho_{s_1},\rho_{s_2}$ are birational self-maps of $X$. Let $A:=\bigcup\limits_{(s_1,s_2)\in R(K)}(\{(s_1,s_2)\}\times A_{(s_1,s_2)}(K))$ which is a dense subset of $S\times S\times X$ by Lemma 2.7(iii).

Now we have that $(s_1,s_2,x)\in\text{Dom}(\rho\circ(m\times\mathrm{id}_{X}))\cap\text{Dom}(\rho\circ(\mathrm{id}_{S}\times\rho))$ for every point $(s_1,s_2,x)\in A$. Moreover, we know $(\rho\circ(m\times\mathrm{id}_{X}))((s_1,s_2,x))=\rho_{m(s_1,s_2)}(x)$ and $(\rho\circ(\mathrm{id}_{S}\times\rho))((s_1,s_2,x))=\rho_{s_1}(\rho_{s_2}(x))$. Both of the notions $\rho_{m(s_1,s_2)}(x)$ and $\rho_{s_1}(\rho_{s_2}(x))$ are well-defined due to our choice of $x$. But we can see that $\rho_{m(s_1,s_2)}=\rho_{s_1}\circ\rho_{s_2}$ just the same as in the proof of Lemma 3.18. So we conclude that $(\rho\circ(m\times\mathrm{id}_{X}))((s_1,s_2,x))=(\rho\circ(\mathrm{id}_{S}\times\rho))((s_1,s_2,x))$ for every $(s_1,s_2,x)\in A$. Therefore we finish the proof by Lemma 2.7(iv).
\end{prf}

\subsection{Step 3}

In this subsection, our goal is to find an open dense subset $S^{(2)}\subseteq S^{(1)}$ such that $\{\rho_{s}|\ s\in S^{(2)}(K)\}$ can be regularized. We say a set $P\subseteq\text{Bir}(X)$ \emph{can be regularized}, if there exist an irreducible projective variety $Y$ and a birational transformation $f:X\dashrightarrow Y$ such that $f\circ P\circ f^{-1}\subseteq\text{Aut}(Y)$.

~

Firstly, we apply the Weil's regularization theorems to the rational group $S$ obtained in subsection 3.2. By Theorem 2.4(i) and Proposition 3.20, we know that there exist a group variety $G_0$ and a birational map $u:G_0\dashrightarrow S$, such that the diagram below is commutative. We \textbf{fix the notation} of $G_0$ and $u$ in the rest of this section.

\[
\xymatrix{G_0\times G_0 \ar[r]^-{m_{G_0}} \ar@{.>}[d]_{u\times u} & G_0 \ar@{.>}[d]^{u} \\ S\times S \ar@{.>}[r]^-{m} & S}
\]

Now the ``group action" $\rho:S\times X\dashrightarrow X$ can induce a dominant rational map $\rho_0=\rho\circ(u\times\text{id}_{X}):G_0\times X\dashrightarrow X$.

\begin{lemma}
The dominant rational map $\rho_0:G_0\times X\dashrightarrow X$ is a rational group action of $G_0$ on $X$ in the sense of Definition 2.2.
\end{lemma}

\begin{prf}
Since the commutativity of the diagram in Definition 2.2(i) is an immediate consequence of Lemma 3.21, we only have to verify that the condition in Definition 2.2(ii) is satisfied. Notice that $\sigma_1:S\times X\dashrightarrow S\times X$ is just $\tilde{\rho}$ (we interchange the $X$ and $S$ in $X\times S$ as before), one can verify that the diagram below is commutative. We remark that in order to be rigorous, one have to verify that the composition of the three maps is $\widetilde{\rho_0}$ directly because $\widetilde{\rho_0}$ may not be \emph{a priori} dominant.
\[
\xymatrix{G_0\times X \ar@{.>}[r]^-{\widetilde{\rho_0}} \ar@{.>}[d]_{u\times\text{id}_{X}} & G_0\times X \\ S\times X \ar@{.>}[r]^-{\sigma_1} & S\times X \ar@{.>}[u]_{u^{-1}\times\text{id}_{X}}}
\]

Now by Proposition 3.8(ii), we conclude that $\widetilde{\rho_0}$ is a birational self-map of $G_0\times X$. Thus we finish the proof.
\end{prf}

We \textbf{fix the notation} of $\rho_0$ in the rest of this section. By Theorem 2.4(ii), we know that there exist an irreducible projective $G_0$-variety $Y$ and a birational map $f:X\dashrightarrow Y$, such that the diagram below is commutative (where $a$ is the group action of $G_0$ on $Y$). We \textbf{fix the notation} of $Y,a$ and $f$ in the rest of this section.
\[
\xymatrix{G_0\times X \ar@{.>}[r]^-{\rho_0} \ar@{.>}[d]_{\mathrm{id}_{G_0}\times f} & X \ar@{.>}[d]^{f} \\ G_0\times Y \ar[r]^-{a} & Y}
\]

The next lemma is an immediate consequence of Lemma 2.10.

\begin{lemma}
The rational group action $\rho_0:G_0\times X\dashrightarrow X$ induces a set of birational self-maps $\{\rho_{0,g_0}'|\ g_0\in G_0(K)\}$ of $X$ (in the way of Lemma 2.10(i)), and this set can be regularized on the birational model $Y$ of $X$.
\end{lemma}

Now we can describe our open dense subset $S^{(2)}\subseteq S^{(1)}$. Since $u:G_0\dashrightarrow S$ is birational, we can find a variety $U$ which is both an open dense subset of $G_0$ and an open dense subset of $S$. Then $u$ can be represented by the open immersions $U\hookrightarrow G_0$ and $U\hookrightarrow S$. We let $S^{(2)}:=U\cap S^{(1)}$ which is an open dense subset of $S$. We \textbf{fix the notation} of $S^{(2)}$ in the rest of this section.

\begin{lemma}
The set $\{\rho_{s}|\ s\in S^{(2)}(K)\}$ can be regularized on the birational model $Y$ of $X$.
\end{lemma}

\begin{prf}
Firstly, we have an ``action" $\rho_U:U\times X\dashrightarrow X$. The rational map $\rho_U$ is dominant and can be gained by two equivalent ways, i.e. $U\times X\hookrightarrow G_0\times X\stackrel{\rho_0}\dashrightarrow X$ and $U\times X\hookrightarrow S\times X\stackrel{\rho}\dashrightarrow X$. Let $s\in S^{(2)}(K)$ be a point, then $s$ lies in $U(K)$. One can verify that both $\rho_{0,s}'$ and $\rho_{s}$ are the self-map induced by $\rho_U$ and $s\in U(K)$. So $\rho_{0,s}'=\rho_{s}$ for every $s\in S^{(2)}(K)$ and hence the result follows from Lemma 3.23. We remark that the notions $\rho_{0,s}'$ and $\rho_{s}$ are well-defined because of Lemma 2.10(i) and the fact that $s\in S^{(1)}(K)$. Moreover, through this note, we always induce self-maps on the fiber from ``actions" in the way of Lemma 2.10(i) and Definition 3.6(ii).
\end{prf}

\subsection{Step 4}

We will finish the proof of Theorem 1.2 in this subsection. In subsection 3.3, we have already seen that there exists an open dense subset $S^{(2)}\subseteq S^{(1)}$ such that $\{\rho_{s}|\ s\in S^{(2)}(K)\}$ can be regularized on the birational model $Y$ of $X$. In this subsection, we will prove that in fact $G\subseteq\text{Bir}(X)$ can be regularized on the birational model $Y$ of $X$. We want to remark that $G$ is a subset of $T(K)$ and may not be a subset of $S(K)$.

The key proposition of this subsection is the result below. Its proof is indeed very similar to the proof of Proposition 3.15(i). In order to be reader-friendly, we will repeat the construction steps as in Step 2. But we will omit some proofs if they are just the same as the corresponding lemmas in subsection 3.2.

\begin{proposition}
Let $g\in G$ be a birational self-map of $X$. Then there exists an open dense subset $U_g\subseteq S^{(1)}$ and a morphism $t_g:U_g\rightarrow S$ such that for each point $s\in U_g(K)$, we have $t_g(s)$ is the closed point in $\bigsqcup\limits_{d=1}^{\infty} C_{\mathrm{dim}(X),d}(X\times X)$ which corresponds to the closed graph of $g\circ\rho_{s}:X\dashrightarrow X$. Notice that $\rho_s$ is well-defined since $s\in S^{(1)}(K)$.
\end{proposition}

~

Now we start our construction steps. We \textbf{fix} an element $g\in G$ in the procedure below.

Let $f_g:=(g\times\text{id}_S)\circ\sigma_1$ be a birational $S$-self-map of $X\times S$. Just the same as in Lemma 3.10, the lemma below can be easily deduced from the fact that $g$ and $\rho_s$ are dominant.

\begin{lemma}
Let $s\in S^{(1)}(K)$ be a closed point. Then the fiber $\mathrm{Dom}(f_g)_{s}$ is an open dense subset of $X$, and the induced rational map $f_{g,s}:X\dashrightarrow X$ equals to $g\circ\rho_{s}$.
\end{lemma}

Now we define $\Gamma\subseteq X\times X\times S$ as the closed graph of the $S$-rational map $f_g$. See the diagram below. We denote $U$ as the open dense subset $\text{Dom}(f_g)\subseteq X\times S$, and the underlines below $X$ are used to indicate that those morphisms are obtained by timing with \emph{that} $\text{id}_{X}$. The closed subvariety $\Gamma$ is the closure of the locally closed subvariety $U\subseteq X\times X\times S$. Notice that the open immersion $U\hookrightarrow\Gamma$ at the top of this diagram makes $U$ into an open dense subset of $\Gamma$.
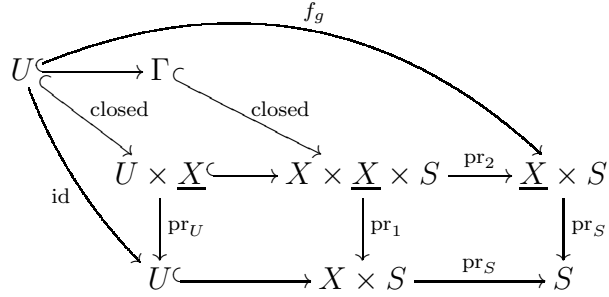
\[
\xymatrix{U \ar@{^{(}->}[r] \ar@/^3pc/[rrrd]^{f_g} \ar@{^{(}->}[rd]^{\text{closed}} \ar@/_/[rdd]_{\text{id}} & \Gamma \ar@{^{(}->}[rd]^{\text{closed}} & & \\ & U\times\underline{X} \ar@{^{(}->}[r] \ar[d]^{\mathrm{pr}_{U}} & X\times\underline{X}\times S \ar[r]^-{\mathrm{pr}_{2}} \ar[d]^{\mathrm{pr}_{1}} & \underline{X}\times S \ar[d]^{\mathrm{pr}_{S}} \\ & U \ar@{^{(}->}[r] & X\times S \ar[r]^-{\mathrm{pr}_{S}} & S}
\]

We denote $p:\Gamma\rightarrow S$ as the surjective structure morphism obtained by the diagram above. Just the same as in Step 2, we can see that there exists an open dense subset $U_g\subseteq S$ such that:
\begin{enumerate}
\item
$U_g$ is a nonsingular open subvariety of $S^{(1)}$.

\item
The morphism $p^{-1}(U_g)\rightarrow U_g$ induced by $p$ is flat.

\item
For each point $s\in U_g$, the fiber $U_{s}\subseteq X_{k(s)}$ can be viewed as an open dense subset of $\Gamma_{s}$.
\end{enumerate}

The next lemma is an immediate consequence of Lemma 3.26 and the part (iii) above.

\begin{lemma} 
\begin{enumerate}
\item
For each closed point $s\in U_g(K)$, the closed subvariety $\Gamma_{s,\mathrm{red}}\subseteq X\times X$ is the closed graph of the birational map $g\circ\rho_{s}:X\dashrightarrow X$.

\item
For each point $s\in U_g$, the fiber $\Gamma_{s}\subseteq(X\times X)_{k(s)}$ is an irreducible closed subscheme of dimension equals to $\mathrm{dim}(X)$ and contains a geometrically reduced open dense subset (which is $U_{s}$).
\end{enumerate}
\end{lemma}

\proof[Proof of Proposition 3.25]
Let $g\in G$ be a birational self-map of $X$ and let $U_g$ be the open dense subset of $S$ obtained in the procedure above. Notice that $U_g\subseteq S^{(1)}$. The proof contains three steps.
\begin{enumerate}
\item
There exists a morphism $t_g':U_g\rightarrow\bigsqcup\limits_{d=1}^{\infty} C_{\mathrm{dim}(X),d}(X\times X)$ such that for each point $s\in U_g(K)$, the closed point $t_g'(s)\in\bigsqcup\limits_{d=1}^{\infty} C_{\mathrm{dim}(X),d}(X\times X)(K)$ corresponds to the closed graph of $g\circ\rho_{s}:X\dashrightarrow X$.

\item
The morphism $t_g':U_g\rightarrow\bigsqcup\limits_{d=1}^{\infty} C_{\mathrm{dim}(X),d}(X\times X)$ factors through $T\subseteq\bigsqcup\limits_{d=1}^{\infty} C_{\mathrm{dim}(X),d}(X\times X)$.

\item
The morphism $t_g':U_g\rightarrow T$ factors through $S$, and thus we obtain the desired morphism $t_g:U_g\rightarrow S$.
\end{enumerate}

The proofs of the three steps above are just the same as the proofs of Lemma 3.12--14. So we omit the detailed proofs. We remark that the key observation of the proof of step (ii) is that $G\cap U_g(K)$ is dense in $U_g$; and the key observation for step (iii) is that $t_g'$ is an injective map on $U_g(K)$.
\endproof

We regard the morphisms $t_g:U_g\rightarrow S$  in Proposition 3.25 as rational self-maps $t_g:S\dashrightarrow S$. Notice that $t_g$ is an injective map on $U_g(K)$, we conclude that the maps $t_g$ are dominant by Lemma 2.7(v). We \textbf{fix the notation} of $U_g,t_g$ as above in the rest of the proof.

~

Now we can finish the proof of Theorem 1.2

\proof[Proof of Theorem 1.2]
We prove that $G\subseteq\text{Bir}(X)$ can be regularized on the birational model $Y$ of $X$. Let $g\in G$ be a birational self-map of $X$. Since $S^{(2)}$ is an open dense subset of $S$ and the self-map $t_g$ is dominant, we know that the inverse image of $S^{(2)}$ under the morphism $t_g:U_g\rightarrow S$ is an open dense subset of $S$. So in particular, there exists a point $s\in(U_g\cap S^{(2)})(K)$ such that $t_g(s)$ also lies in $S^{(2)}(K)$. Now the closed point $t_g(s)\in\bigsqcup\limits_{d=1}^{\infty} C_{\mathrm{dim}(X),d}(X\times X)(K)$ corresponds to the closed graph of both $\rho_{t_g(s)}$ and $g\circ\rho_s$. So we have $\rho_{t_g(s)}=g\circ\rho_s$. Therefore we finish the proof by Lemma 3.24.
\endproof

\section{Another approach using Hilbert schemes}

In this section, we describe a proof of Theorem 1.2 using Hilbert scheme as the moduli space. Since the main procedure is just the same, we will only describe the modifications and complements in each step.

Firstly, we recall the concept of Hilbert schemes. The original source is Grothendieck's FGA and one can consult $\cite[\text{Chapter}\ 5]{FGA05}$ for a reference.

\begin{theorem}(Hilbert scheme)

Let $X$ be a projective variety and let $L$ be a very ample line bundle on $X$. Fix a polynomial $P\in\mathbb{Q}[\lambda]$. Then the functor from $(\mathrm{Varieties)}^{op}$ to $\mathrm{(Sets)}$ as below can be represented by a projective variety $\mathrm{Hilb}_{P,L}(X)$ (we consider the reduced structure of the usual Hilbert scheme because we only deal with varieties in this note):
$$
S \mapsto \left\{
\begin{array}{@{}l@{}}
  F\subseteq X\times S\text{ is a closed}\\
  \text{subscheme of }X\times S\\
\end{array}
\;\left|\;\;
\begin{array}{@{}l@{}}
  \text{$\forall s\in S$, the Hilbert polynomial of $F_s\subseteq X_{k(s)}$}\\
  \text{with respect to $L$ is $P$; and the morphism}\\
  \text{$F\rightarrow S$ is flat.}\\
\end{array}
\right.
\right\}
$$
\end{theorem}

In order to carry out the argument in Section 3 for Hilbert schemes, we need a result to guarantee that the closed graphs of the elements in $G$ lie in finitely components of the Hilbert scheme. This result is stated in the proposition below and it is an immediate consequence of $\cite{Starr}$.

\begin{proposition}
Let $X$ be a projective variety and let $L$ be a very ample line bundle on $X$. Then for a fixed positive integer $d$, there are only finitely many possibilities for the Hilbert polynomial $h_Y$ of an irreducible closed subvariety $Y\subseteq X$ whose degree $\mathrm{deg}(Y)=d$. Both of the degree and the Hilbert polynomial are calculated by $L$.
\end{proposition}

Now we start to modify the argument in Section 3. We modify Step 0 now, and the modifications for Step 1--4 will take place in subsections 4.1--4. We remark that basically we will inherit \textbf{every} notation in Section 3. The notations without modification will remain the same (or more precisely, be constructed in the same way) as in Section 3 and the notations with modification will also be very similar and play the same role as in Section 3. Also, we will only describe the changes at the points where changes are needed. The arguments that we do not mention should still go well.

~

In this section, we let $T$ be the closure of $\{\Gamma_f|\ f\in G\}\subseteq\bigsqcup\limits_{P\in\mathbb{Q}[\lambda]}\mathrm{Hilb}_{P,L_0}(X\times X)(K)$ in $\bigsqcup\limits_{P\in\mathbb{Q}[\lambda]}\mathrm{Hilb}_{P,L_0}(X\times X)$. Then Proposition 4.2 guarantees that $T$ is indeed a projective variety because $G$ is of bounded-degree. We still let $S$ be the maximal dimension part of $T$ as in Section 3. Now Proposition 3.1 should read as follows.

\begin{proposition}
We restrict the universal family on $\bigsqcup\limits_{P\in\mathbb{Q}[\lambda]}\mathrm{Hilb}_{P,L_0}(X\times X)$ to $S$ and get a closed subscheme $F\subseteq X\times X\times S$. Then $F$ is flat over $S$ and the following holds:
\begin{enumerate}
\item
Let $s\in S(K)$ be a closed point. Then the fiber $F_s\subseteq X\times X$ is tautologically the closed subscheme of $X\times X$ which corresponds to the closed point $s\in\bigsqcup\limits_{P\in\mathbb{Q}[\lambda]}\mathrm{Hilb}_{P,L_0}(X\times X)(K)$.

\item
If $g\in G\cap S(K)$, then $F_g=\Gamma_g$ as a closed subvariety of $X\times X$.
\end{enumerate}
\end{proposition}

We also mention that this flat family over $S$ has good properties.

\begin{lemma}
\begin{enumerate}
\item
$\{s\in S|\ F_s\text{ is geometrically integral}\}$ is an open dense subset of $S$.

\item
$F$ is indeed a closed sub\emph{variety} of $X\times X\times S$, i.e. it is reduced.
\end{enumerate}
\end{lemma}

\begin{prf}
\begin{enumerate}
\item
This set is open in $S$ because of $\cite[\text{Th\'eor\`eme}\ 12.2.4\text{(viii)}]{GD66}$ and it is dense in $S$ because it contains $G\cap S(K)$.

\item
Since $F$ is flat over $S$, this follows from $\cite[4.3.1,\text{Proposition}\ 3.8]{Liu02}$ and part (i).
\end{enumerate}
\end{prf}

\subsection{Changes in Step 1}

There is only one minor change for the arguments in this step. In view of Lemma 4.4(i), we may just let $S^{(1)}$ be the open dense subset $\{s\in S|\ F_s\text{ is geometrically integral}\}\cap\text{pr}_{S}(S_1)$ of $S$. And other things just remain the same.

\subsection{Changes in Step 2}

Recall that the main goal of this step is to construct the rational group law on $S$ and our idea is to construct an appropriate family $\Gamma$ over $S\times S$ in order to induce the rational group law. Since the functor of Hilbert schemes is scheme-theoretic in contrast to the functor of open Chow varieties, a little more work is needed in order to show that $\Gamma$ is also capable for its responsibility in the approach of this section. We will modify the prior ``domain of definitions" $Q,Q_1,Q_2$ of the rational maps $m,m_1,m_2:S\times S\dashrightarrow S$ here. Now we show the modification procedure for $Q$ as an example. One can modify $Q_1$ and $Q_2$ in exactly the same way.

\begin{lemma}
For every generic point $\eta\in S\times S$, the fiber $\Gamma_{\eta}$ is geometrically reduced.
\end{lemma}

\begin{prf}
Firstly, the generic fiber $\Gamma_{\eta}$ is reduced as both $\Gamma$ and $S\times S$ are varieties. Notice that by the proof of Lemma 2.14, we see that $U_{\eta}\subseteq\Gamma_{\eta}$ is an open dense geometrically reduced subset. Then we get the conclusion by $\cite[\text{Lemma}\ 33.6.8\text{(3) (04KS)}]{Stacks}$.
\end{prf}

Now we can modify the prior ``domain of definition" $Q$ of the rational group law $m:S\times S\dashrightarrow S$ as follows. Combining the generic flatness theorem, Lemma 2.12 and Lemma 2.14 and $\cite[\text{Th\'eor\`eme}\ 12.2.4\text{(v)}]{GD66}$ and Lemma 4.5, we can choose an open dense subset $Q\subseteq S\times S$ such that:
\begin{enumerate}
\item
$Q\subseteq S^{(1)}\times S^{(1)}$.

\item
The morphism $p^{-1}(Q)\rightarrow Q$ induced by the surjective structure morphism $p:\Gamma\rightarrow S\times S$ is flat.

\item
For each point $q\in Q$, the fiber $U_q\subseteq X_{k(q)}$ can be viewed as an open dense subset of $\Gamma_q$; and the fiber $\Gamma_q$ is geometrically reduced.
\end{enumerate}

Then we have that for each closed point $q=(s_1,s_2)\in Q(K)$, the fiber $\Gamma_q\subseteq X\times X$ is the closed graph of the birational map $\rho_{s_1}\circ\rho_{s_2}:X\dashrightarrow X$. The reader can modify $Q_1$ and $Q_2$ in the way same as above, and then the remaining arguments in Step 2 go well. One can just literally change the open Chow varieties in the argument into Hilbert schemes.

\subsection{Changes in Step 3}

The arguments in this step remain valid line by line and no changes is needed to make.

\subsection{Changes in Step 4}

The main part of the argument in Step 4 is just a repeat of Step 2, so the changes are also exactly the same as in subsection 4.2. One can prove an analogue of Lemma 4.5 and modify the prior ``domain of definition" $U_g$ in the way as in subsection 4.2. We shall leave this little work to the reader. And then we can finish this proof of Theorem 1.2 use the same argument as in the end of Section 3.

\section{Bounded-degree self-maps in arithmetic dynamics}

As one can imagine, the bounded-degree self-maps are the ``simplest" self-maps of an irreducible projective variety. In differential dynamics, one measures the complexity of a self-map by the entropy. In arithmetic dynamics, there is a fundamental concept called \emph{dynamic degree} which is an analogue of the entropy in differential dynamics. There are many references for this concept in the literature, see for example $\cite{Dang20},\cite{Tru20}$ and $\cite[\text{Section}\ 2.1]{Xie23}$. We will introduce the concept of the \emph{degree sequence} and the definition of bounded-degree self-maps in this section. Then we deduce Corollary 1.3 from Theorem 1.2.

~

Through this section, we let $X$ be an irreducible projective variety and let $f$ be a dominant rational self-map of $X$. Let $L\in\text{Pic}(X)$ be a big and nef line bundle. We consider the graph $\Gamma_f\subseteq X\times X$ which is an irreducible closed subvariety. Let $\pi_1,\pi_2:\Gamma_f\rightarrow X$ be the two projections. Then $\pi_1$ is a birational proper morphism and $\pi_2$ is surjective proper. We define the \emph{first degree} $\text{deg}_{1,L}(f)$ of $f$ with respect to $L$ as the intersection number $(\pi_2^{*}(L)\cdot\pi_1^{*}(L)^{\text{dim}(X)-1})$ on $\Gamma_f$. Then we get a sequence $\{\text{deg}_{1,L}(f^n)|\ n\in\mathbb{N}\}$ of positive integers.

For two sequences $\{a_{n}\},\{b_{n}\}\in(\mathbb{R}_{\geq1})^{\mathbb{N}}$, we say that they have the same speed of growth if $\{\frac{a_{n}}{b_{n}}|\ n\in\mathbb{N}\}$ has an upper bound and a positive lower bound. Let $\mathrm{deg}_{1}(f)$ be the class of the speed of growth of the sequence $\{\mathrm{deg}_{1,L}(f^{n})|\ n\in\mathbb{N}\}$, which by $\cite[\text{Theorem}\ 1\text{(ii)}]{Dang20}$ is irrelevant with the choice of the big and nef line bundle $L$. Notice that although $\cite[\text{Theorem}\ 1\text{(ii)}]{Dang20}$ was stated for normal projective variety $X$, the result also holds for arbitrary irreducible projective variety $X$ because one can pass to the normalization of $X$. Then we can abuse notation and say that $\mathrm{deg}_{1}(f)$ is the \emph{degree sequence} of $f$. We also remark that in fact $\mathrm{deg}_{1}(f)$ remains the same on different birational models. See $\cite[\text{top of p. 1269}]{Dang20}$.

Now we state the definition of bounded-degree self-maps.

\begin{definition}
Let $X$ be an irreducible projective variety. We say a dominant rational self-map $f:X\dashrightarrow X$ is of bounded-degree if $\mathrm{deg}_{1}(f)$ is a bounded sequence. 
\end{definition}

The main proposition of this section is as follows.

\begin{proposition}
Let $X$ be an irreducible projective variety and let $f$ be a bounded-degree self-map of $X$. Then $f$ is a birational self-map and $\{f^{n}|\ n\in\mathbb{Z}\}\subseteq\mathrm{Bir}(X)$ is a bounded-degree subgroup in the sense of Definition 1.1.
\end{proposition}

We start with a technical lemma in the intersection theory.

\begin{lemma}
Let $X$ be an irreducible projective variety of dimension $n\geq2$ and let $L_1,\dots,L_n$ be nef line bundles on $X$. Then we have $(L_1\cdot L_1\cdot L_3\cdots L_n)\cdot(L_2\cdot L_2\cdot L_3\cdots L_n)\leq(L_1\cdot L_2\cdot L_3\cdots L_n)^2$.
\end{lemma}

\begin{prf}
Since nef line bundles are the limit of ample line bundles, we may assume that $L_1,\dots,L_n$ are ample line bundles without loss of generality. We prove by induction on $n=\text{dim}(X)$. When $n=2$, one can resolve the singularities of $X$ and the result follows from the Hodge index theorem.

For the general case, one can assume that $L_n$ is very ample without loss of generality and complete the reduction procedure by a Bertini-type result $\cite[\text{Lemma}\ 37.32.3\ \text{(0G4F)}]{Stacks}$. We need this Bertini-type result in order to guarantee that we always face to irreducible projective varieties during the reduction procedure.
\end{prf}

Now we can prove Proposition 5.2.

\proof[Proof of Proposition 5.2]
Let $L\in\text{Pic}(X)$ be an ample line bundle. The bounded-degree condition of $f$ says that $\{(\pi_{2,n}^{*}(L)\cdot\pi_{1,n}^{*}(L)^{\text{dim}(X)-1})|\ n\in\mathbb{N}\}\subseteq\mathbb{Z}_+$ is bounded in which $\pi_{1,n},\pi_{2,n}:\Gamma_{f^n}\rightarrow X$ are the two projections. Recall that we defined the notion $\text{deg}_{L}(f)$ in Introduction before Definition 1.1. We calculate that $\text{deg}_{L}(f^n)=(\pi_{2,n}^{*}(L)+\pi_{1,n}^{*}(L))^{\text{dim}(X)}$ for every nonnegative integer $n$. Notice that since $\pi_{1,n}$ are birational morphisms, we have that $\pi_{1,n}^{*}(L)^{\text{dim}(X)}=L^{\text{dim}(X)}$ is a constant irrelevant to $n$. So by Lemma 5.3, we can firstly conclude that $\{\text{deg}_{L}(f^n)|\ n\in\mathbb{N}\}\subseteq\mathbb{Z}_+$ is a bounded set.

Now in particular, we can see that the set $\{\pi_{2,n}^{*}(L)^{\text{dim}(X)}|\ n\in\mathbb{N}\}$ is bounded. But since $\pi_{2,n}$ is a surjective morphism between irreducible projective varieties of degree equals to $\text{deg}(f)^n$, we conclude that $\pi_{2,n}^{*}(L)^{\text{dim}(X)}=\text{deg}(f)^n\cdot(L^{\text{dim}(X)})$ and hence we must have $\text{deg}(f)=1$. Thus we conclude that $f$ is a birational self-map of $X$.

Now we only have to prove that $\{\text{deg}_{L}(f^n)|\ n\in\mathbb{Z}\}$ is a bounded set. Notice that for every positive integer $n$, the closed subvarieties $\Gamma_{f^{n}},\Gamma_{f^{-n}}\subseteq X\times X$ are ``symmetric" in the sense that they correspond to each other under the swapping automorphism $X\times X\stackrel{\sim}\rightarrow X\times X$. Therefore by the definition of the notion $\text{deg}_{L}(f)$, we know that $\text{deg}_{L}(f^n)=\text{deg}_{L}(f^{-n})$ for every positive integer $n$. Hence we finish the proof as we have seen that $\{\text{deg}_{L}(f^n)|\ n\in\mathbb{N}\}$ is bounded.
\endproof

Lastly, we prove Corollary 1.3.

\proof[Proof of Corollary 1.3]
This corollary follows immediately from Theorem 1.2 and Proposition 5.2. Notice that by the argument in $\cite[\text{top of p. 1269}]{Dang20}$, we know that the degree sequence $\text{deg}_{1}(f)$ remains the same after changing of birational models.
\endproof

\bibliographystyle{alpha}
\bibliography{reference}

\end{spacing}
\end{document}